Andrei Rodin (rodin@ens.fr)

Ecole Normale Supérieure

# Identity and Categorification

**Content:**

# 1. Paradoxes of Identity and Mathematics

Changing objects (of any nature) pose a difficulty for the metaphysically-minded logician known as the *Paradox of Change* . Suppose a green apple becomes red. If *A* denotes the apple when green, and *B* when it is red then *A=B* (it is the same thing) but the properties of *A* and *B* are different: they have a different color. This is at odds with the *Indiscernibility of Identicals* thesis according to which identical things have identical properties. A radical solution - to explain away and/or dispense with the notion of change altogether was first proposed by Zeno around 500 BC and remains popular among contemporary philosophers (who often appeal to the relativistic spacetime to justify Eleatic arguments). Unlike physics, mathematics appeared to provide support for the Eleatic position: for some reason people were more readily brought to accept the idea that mathematical objects did not change than to accept a similar claim about physical objects - in spite of the fact that mathematicians had always talked about variations, motions, transformations, operations and other process-like notions just as much as physicists.

The Paradox of Change is the common ancestor of a family of paradoxes of identity which might be called *temporal* because all of them involve objects changing in time[1]. However *time* is not the only cause of troubles about identity: *space* is another. The *Identity of Indiscernibles* (the thesis dual to that of the *Indiscernibility of Identicals*) says that perfectly like things are identical. According to legend in order to demonstrate this latter thesis, Leibniz challenged a friend during a walk to find a counter-example among the leaves of a tree. Although there are apparently no perfect doubles among material objects, mathematics appears to provide clear instances immediately: think about two (different) points. But the example of geometrical space brings another problem: either the *Identity of Indiscernibles* thesis is false or our idea of perfect doubles like points is incoherent. In what follows I shall refer to this latter problem as the *Paradox of Doubles*. Mathematics looks more susceptible to this paradox than physics. However nowadays mathematics and physics are so closely entwined it is hardly possible to isolate difficulties in one

---

[1] *Chrisippus' Paradox*, *Stature*, *The Ship of Theseus* belong to this family. See Deutsch (2002)

discipline from those in the other. Were she living today, Leibniz's friend might meet his challenge by mentioning the indiscernibility of particles in Quantum Physics [2].

Category theory provides an original understanding of identity in mathematics, which takes seriously the idea that mathematical objects are, generally speaking, *variable* and handles the problem of doubles in a novel way. Category theory does not *resolve* paradoxes of identity of the above form; rather it provides a setting where paradoxes in such a form do not arise. The Category-theoretic understanding of identity in mathematics may have important consequences in today's mathematically-laden physics and hence (assuming some form of scientific realism) for the notion of identity in a completely general philosophical setting. In this paper I explore this new understanding of identity in Category Theory, leaving its implications outside mathematics for a future study.

The paper is organized as follows. First, I consider some difficulties about the notion of identity in mathematics, providing details and examples. Then I briefly review some attempts to overcome these difficulties. I pay particular attention to the account of identity in mathematics proposed by Frege and afterward developed by Russell, which remains standard in the eyes of many philosophers. Then I consider an alternate approach to identity in mathematics, which dates back to Greek geometry but made a new appearance in 19th century and later developed in Category theory. I consider the issue of identity in Category theory starting with general remarks and then coming to more specific questions concerning fibred categories and higher categories. Finally I suggest a way of thinking about categories, which implies deversification of the notion of identity and revision of Frege's assumption that identities must be fixed from the outset.

## 2. Mathematical Doubles

The example of two distinct points *A*,*B* (Fig.1) does not, it is usually argued, refute the *Identity of Indiscernibles* because the two points have different *relational* properties: in Fig. 1, *A* lies to the left of *B* but *B* does not lie to the left of itself[3]:

---

[2] French (1988)

*A* •                               • *B*

Fig.1

(The difference in the relational properties of *A* and *B* amounts to saying that the two points have *different positions*.)

However the example can be easily modified to meet the argument. Consider two *coincident* points (Fig.2): now *A* and *B* have the *same* position.

*A=B*

•

Fig.2

It might be argued that coincident points are an exotic case, one which can and should be excluded from mathematics via its logical regimentation. But this is far from evident - at least if we are talking about classical Euclidean geometry. For one of the basic concepts of Euclidean geometry is *congruence,* and this notion (classically understood) presumes coincidence of points: figures *F*, *G* are *congruent* iff by moving *G* (without changing its shape and its size) one can make *F* and *G* *coincide* point by point.

The fact that geometrical objects may coincide differentiates them significantly from material solids like chairs or Democritean atoms. The supposed *impenetrability* of material solids counts essentially in providing their identity conditions (Lucas 1973). Thus, identity works differently for material atoms and geometrical points[4].

---

[3] These relational properties of the two points depend on their shared space: the argument doesn't go through for points living on circle. I owe this remark to John Stachel.

[4] This fact shows that Euclidean geometrical space cannot be viewed as a realistic model of the space of everyday experience as is often assumed. One needs the third dimension of physical

We see that the alleged contradiction with the *Identity of Indiscernibles* is not the only difficulty involved here. Indeed the whole question of identity of points becomes unclear insofar as they are allowed to coincide. Looking at Fig.2 we have a surprising freedom in interpreting "=" sign. Reading "=" as identity we assume that *A* and *B* are two different names for the same thing. Otherwise we may read "=" as specifying a coincidence relation between the (different) points *A* and *B*. It is up to us to decide whether we have only one point here or a family of superposed points. The choice apparently has little or no mathematical sense. One may confuse coincidence with identity here without any risk of error in proofs. However this does not mean that one can just assimilate the notions of identity and coincidence. For identity so conceived would be very ill-behaved, allowing for the merger of different things into one and the splitting of one into many. (Consider the fact that Euclidean space allows for the coincidence of *any* point with any other through a suitable motion.) Perhaps it would be more natural to say instead that the relations of coincidence and identity while not identical in general, coincide in this context?

For an example from another branch of elementary mathematics consider this equation: 3=3. Just as in the previous case there are different possible interpretations of the sign "=" here. One may read "=" either as identity, assuming that 3 is a unique object, or as a specific relation of equality which holds between different "doubles" (copies) of 3. Which option is preferable depends on a given context. There is a unique natural number $x$ such that $2<x<4$; $x=3$. Here "=" stands for identity. But when one thinks about the sum 3+3 or the sequence 3,3,3,... it is convenient to think of the 3s as many. In this latter case 3=3 still holds but now "=" is being read as equality rather

---

space to establish in practice the relation of congruence between (quasi-) 1- and 2-dimensional material objects through the application of a measuring rod or its equivalent.

[6] The unit of a group *G* is defined as the element $1 \in G$ such that for any element $x \in G$ (including 1 itself) $1 \otimes x = x \otimes 1 = x$, where $\otimes$ is the group operation. The existence of 1 is guaranteed by definition but its uniqueness is proved. Suppose 1' is another element of the group satisfying the same condition: $1' \otimes x = x \otimes 1' = x$. Then taking first $x=1$, and then $x=1'$ we have $1' \otimes 1 = 1 \otimes 1' = 1 = 1'$. This argument justifies the use of the definite article in the expression "*the* unit of *G*".

than identity. Again the choice looks like a matter of convenience rather than of theoretical importance.

Similarly, in one sense *cube* is a particular geometrical object, while in a different sense there exist (in some suitable sense of "exist") many cubes. When one proves that there exist exactly 5 different regular polyhedrons, and says that the *cube* is one of them, one speaks about the cube in the first sense. When one considers a geometrical construction, which comprises several cubes, one thinks about the cubes in the second sense. However no distinction between the two meanings of the term "cube" can be found in standard textbooks, and it is not even clear whether such distinction can be sharply made.

The above examples might make one think that the notion of identity simply plays no significant role in mathematics. 2x2=4 remains true independently of whether the sign "=" is read as equality, or as identity, whether equality is treated as identity, or identity is weakened to equality. It looks as if here one may choose one's interpretation according to personal taste or preferred philosophical position.

However such a liberal attitude to identity in mathematics looks suspicious from the logical point of view. Claims of existence and uniqueness of mathematical objects satisfying given descriptions (definitions) play an important role in mathematics. Such a claim means that a given description indeed picks out (identifies) an object, not just a property. The standard definition of the *unit* of a given group *G* (also often called the *identity of G*) is an example.[6] Obviously a claim that such-and-such an object is unique makes sense only if its identity conditions are fixed. But as we have seen they may in fact be very loose. It is clear that 3 is the only natural number bigger than 2 and smaller than 4 but it is not clear that 3 indeed refers to an unique object. But how can mathematics hang together as a body of knowledge if it apparently does not meet Quine's "no entity without identity" requirement?

There are several ways to approach this problem. I now explore them.

## 3. Types and Tokens

The remedy, which readily comes to mind on the part of anyone familiar with contemporary Analytic metaphysics, is that of the type/token distinction. Consider another example, which prima facie looks very like the above mathematical cases. There are 26 letters in the English alphabet, and the letter *a* is one of them. In the last phrase the letter *a* is referred to as a particular thing, namely a particular letter of the alphabet. But in this phrase itself there are five such things. Hence the letter *a* is not a particular thing. The standard way of dissolving this puzzle is to say that here we have one *a*-type and five *a*-tokens.

In explaining the distinction, one starts from tokens: an a-token is a piece of paper with typographic pigment, or another material object (e.g. a piece of printer's type) representing the letter *a*. Obviously *a*-tokens are many. The second step is to explain what the *a*-type is. Intuitively it is what all and only *a*-tokens share in common (typically a certain shape). To explain the notion of type better than this is not an easy task, and it involves old and hard metaphysical questions as well as complicated logical problems, which I shall not enter into here. Let me show instead that the type/token distinction doesn't fix the problem of identity of mathematics anyway: whatever mathematical types might be they do not correspond to well-distinguishable tokens.

*The* natural number **3** (which I write in bold for further references) indeed looks like a type but the 3s, which we find in the series 3,3,.. or in the formula 3+3 do not look like tokens from the viewpoint of standard examples (like particular chairs). For formula 3+3 may be applied to many different situations: one might add 3 chairs to 3 chairs, 3 points to 3 points, or even (taking a liberal attitude) 3 chairs to 3 points[7]. Arguably such application amounts to instantiation of both 3s (in formula 3+3) by certain sets of objects. That is certainly not how good tokens behave: the fact that types can be instantiated but tokens cannot is essential; if we allow for the instantiation of tokens by other tokens we either lose the type/token distinction or must provide it with a new relational sense (which looks like an interesting project but I cannot pursue it here).

---

[7] The last example shows that the typification certainly matters here but this does not change the argument.

The case of points (or more structured geometrical figures like triangles) at first sight looks more promising. Apparently points are well-distinguishable tokens of the same type. Unlike the case of numbers it is common in mathematics to denote different point-tokens by different labels such as *A* and *B*. However this works only until coincident points are taken into consideration. For in the case of coincident points we cannot distinguish a singular point-token from a "stock" of point-tokens. It is tempting in this case to think of the stock of points as a "place" occupied by a family of singular point-tokens. But this again involves a reiteration of the type/token distinction on another level as in the case of 3-tokens. Point-locations initially considered as tokens can themselves be instantiated by second-order tokens stocked there. Once again this destroys the usual distinction between point-tokens and the point-type. It is a condition of acting as a (classical) token that the object so acting have determinate identity conditions - as concrete symbols like printed numerals do. But our hypothetical number- and point-tokens do not meet this condition. So the type/token distinction (at least in its usual form) does not help us to handle the identity issue in mathematics. (This also makes me doubt how well it works outside mathematics.)

## 4. Frege and Russell on The Identity of Natural Numbers

Frege considered it a principal task of his logical reform of arithmetic to provide absolutely determinate identity conditions for the objects of that science, i.e. for numbers. Referring to the contemporary situation in this discipline he writes in the Introduction to his (1884):

How I propose to improve upon it can be no more than indicated in the present work. With numbers ... it is a matter of fixing the sense of an identity.   (English trans. 1953, p.Xe)

Frege makes the following critically important assumption : identity is a general logical concept, which is not specific to mathematics.  In (1884) Frege says:

It is not only among numbers that the relationship of identity is found. From which it seems to follow that we ought not to define it specially for the case of numbers. We should expect the concept of identity to have been fixed first, and that then from it together with the concept of number it must be possible to deduce when numbers are identical with one another, without there being need for this purpose of a special definition of numerical identity as well. (Trans. 1953, p.74e)

In a different place Frege says clearly that this concept of identity is absolutely stable across all possible domains and contexts:

Identity is a relation given to us in such a specific form that it is inconceivable that various forms of it should occur[8] (Frege 1903; edition 1962, p.254)

Frege's definition of natural number, as modified in Russell (1903) later became standard[9]. I present it here informally in Russell's simplified version. Intuitively the number 3 is what all collections consisting of three members (trios) share in common. Now instead of looking for a common form, essence or type of trios let us simply consider all such things together. According to Frege and Russell the collection (class, set[10]) of all trios *just is the number* 3. Similarly for other numbers.

Isn't this construction circular? Frege and Russell provide the following argument which they claim allows us to avoid circularity here: given two different collections we may learn whether or not they have the same number of members without knowing this number and even without the notion of number itself. It is sufficient to find a one-one correspondence between members of two given collections. If there is such a correspondence, the two collections comprise the same number of members, or to avoid any reference to numbers we can say that the two collections are *equivalent*. I shall follow current usage in calling this equivalence *Humean*.[11] Now we check that this relation is indeed an equivalence in the usual sense, and define natural numbers as equivalence classes under this relation.

---

[8] "Die Identitaet ist eine so bestimmt gegebene Beziehung, dass nicht abzusehen ist, wie bei ihr verschiedene Arten vorkommen koennen."

[9] See, for example Fraenkel (1966), p 10.

[10] Following Russell (1903) I use here words *class*, *collection*, and *set* interchangeably ignoring their technical meanings if any. This terminological freedom is helpful for rethinking the concept of set (or class etc.) without smuggling in ready-made solutions through the existing terminology.

[11] Hume (1978), book 1, part 3, sect. 1

This definition reduces the question of identity of numbers to that of identity of classes. This latter question is settled through the axiomatization of set theory in a logical calculus with identity. Thus Frege's project is realized: it has been seen how the logical concept of identity applies to numbers. (In fact this doesn't work that smoothly as I show in the section **10** below.) In an axiomatic setting "identities" in Quine's sense (that is, identity conditions) of mathematical objects are provided by an axiom schema of the form

$\forall x \forall y (x=y \leftrightarrow \_\_\_)$,

called in Keranen (2001) the *Identity Schema* (IS)[12]. This does not resolve the identity problem though because any given system of axioms, generally speaking, has multiple models[13]. The case of isomorphic models is similar to that of equal numbers or coincident points (naively construed): there are good reasons to think of isomorphic models as one and there is also good reason to think of them as many. So the paradox of mathematical "doubles" reappears. It is a highly non-trivial fact that different models of Peano arithmetic, ZF, and other important axiomatic systems are not necessarily isomorphic. Thus logical analysis à la Frege-Russell certainly clarifies the mathematical concepts involved but it does not settle the identity issue as Frege believed it did. In the recent philosophy of mathematics literature the problem of the identity of mathematical objects is usually considered in the logical setting just mentioned: either as the problem of the non-uniqueness of the models of a given axiomatic system or as the problem of how to fill in the Identity Schema. For my present purposes it is important, however, to return to the problem in its original "informal" version, which inspired Frege and Russell 100 years ago. Such a return to the starting point is, in my view, helpful and perhaps necessary if one wishes (as I do) to consider the Category-theoretic approach to identity discussed in this paper as a viable alternative to the approach taken by Frege, Russell and their followers. At the first glance the Frege-Russell proposal concerning the identity issue in mathematics seems judicious and innocent

---

[12] See section **10** below for an example.

[13] Benacerraf (1965)

(and it certainly does not depend upon the rest of their logicist project): to stick to a certain logical discipline in speaking about identity (everywhere and in particular in mathematics). The following historical remark shows that this proposal is not so innocent as it might seem.

**5. Plato**

Given a sequence like 3,3,3... mathematicians conveniently talk about multiple "copies" of the same number (similarly about copies of a given set, or space) . Such talk about "copies" carries echoes from Plato. A glance at Plato's philosophy of mathematics[14] shows some features which might be attractive for a mathematician resistant to the logical regimentation of talk of identity in different contexts proposed by Frege and Russell. If I understand Plato correctly, according to him identity applies only to the immutable *ideas*, and only ideas *exist*. (So Plato's view in this respect is in accord with Quine's dictum about "no entity without identity".) Material things don't exist but *become*   ( they change, come into and go out of being ) and hence have no proper identities: this is another possible way out of the Paradox of Change. Mathematical things occupy an intermediate position between material stuff and ideas: they involve a weaker sort of becoming and a softer form of identity. In the case of numbers such "soft identity" is equality. Things in the three layers of Plato's ontology are partially ordered by "distorted copying" where ideas are the maximal elements,  mathematical objects are distorted copies of ideas, and material objects are distorted copies of mathematical objects (and hence also of ideas). The distortion of self-identical ideal numbers amounts to their replacement by families of equal mathematical numbers. For example, there is a unique ideal number **3** and an indefinite number of its equal mathematical copies. In other words numbers in mathematics are defined up to equality but not up to identity.

I cannot provide here full justification of this reading of Plato and give only the following hints referring to my (2003a) for a systematic treatment. There are multiple passages where Plato speaks of "*X* itself", "*X* (thought of) through itself" ($\kappa\alpha\vartheta$@$\upsilon\tau\omega$) and "Idea of *X*" interchangeably or explains the latter through the former. For example in *Symposium* (210-211) Plato does this with the notion of Beauty, and in *Phedon* (96-103) with number 2. (In this latter dialog Socrates rejects the view that 2 could be thought of as sum of two units pointing to the fact that 2 can be equally obtained through division of given unit into two halves. Since each of the two operations is the reverse of the other none of them can be viewed as bringing 2 about. So one needs to think of the idea of 2 independently of operations of this sort.) I interpret these passages in the sense, which seems me straightforward: "identity to itself" applies to ideas but neither to material things, nor to mathematical things (as they are usually thought of). To see that Plato's "idea of 2" is indeed something else than mathematical number see last chapters of Aristotle's *Metaphysics* where the author criticizes the "Unwritten doctrine" developed by Plato in the later period of his life (Findlay 1974). Here the distinction between ideal and mathematical numbers is made explicit. Aristotle stresses the fact that each ideal number is unique while their mathematical copies are many (*Met*. 987b) and the fact that ideal numbers are not a subject of arithmetical operations (*Met*. 1081a-1082b).

---

[14] Not to be confused with "Mathematical Platonism" in the sense of Balaguer (1998) and many others, which has little if anything to do with historical Plato. For an introduction to Plato's philosophy of mathematics see Pritchard (1995).

Thus unlike Frege Plato does not suppose that the notion of identity applies to whatever there is (or whatever occurs) indiscriminately. Instead Plato thinks of identity as a specific property of things he calls *ideas* and notices the fact that in mathematics the identity requirement is relaxed. In what follows I shall show that this Platonic insight is particularly appealing in the context of our contemporary Category-theoretic mathematics.

### 6. Equality and Equivalences in Geometry

Plato hints at the following division of labor: mathematicians work on equalities whilst philosophers take care of identities. In the case of arithmetic this is exactly what mathematicians (and philosophers like Frege) have been doing for centuries. In geometry however the situation is more complicated because *equality* in this discipline may mean - and historically did mean - different things.

Euclid uses the term "equality" ($\iota\sigma o\nu$) in the sense of equicompositionality (of plane geometrical figures ) but there are other equivalencies in geometry, which may be considered as better "working substitutes for identity": for example congruence, (geometric) similarity, and affinity. For there is a sense in which the "same figure" means a figure of the same shape and the same size, and there is another sense in which it means only a figure of the same shape, and the notion of "same shape" can itself also be specified in different ways. In addition geometry unlike arithmetic allows for the identification of its objects (of geometrical figures) by directly naming them, usually through naming of their most important points. This allows us to distinguish two different triangles *ABC* and *A'B'C'* which are the "same" in any of above senses. There is apparently no clear argument, which would allow us to choose one of these senses of "the same" as basic and eliminate the others as an abuse of the language. In particular, as I have shown in section **3** above, pointwise naming of figures cannot do this job. So the situation in geometry (even classical geometry!) is exactly that which Manin (2002) describes for a different purpose: There is no equality in mathematical objects, only equivalences.

## 7. Definitions by Abstraction

To pursue his project of reducing the various informal meanings of "the same" in mathematics to a standard notion of identity captured in a universal logic Frege proposed the method of "definition by abstraction". In his (1884) Frege gives the following example of such definition:

The judgment "line *a* is parallel to line *b*", or, using symbols *a//b* , *can be taken as identity*. If we do this, we obtain the concept of direction, and say: "the direction of line *a* is identical with the direction of line *b*". Thus we replace the symbol // by the more generic symbol =, through removing what is specific in the content of the former and dividing it between *a* and *b*. (reprint 1964, p. 74$^e$, italic mine)

Notice that the procedure as described here by Frege involves a change of notation: in the formula *a=b* the symbols *a,b* no longer stand for lines but denote the same direction. Calling this formal procedure *definition by abstraction* Frege suggests its interpretation. The idea is that the procedure picks out a property common to all members of a given equivalence class. In section **11** I shall show that this procedure can be interpreted differently.

As our earlier quotations from Frege (1884) clearly show, in treating an equivalence *E* "as identity" Frege does not mean to replace identity by something else. He aims at the exact opposite: to introduce identity where mathematicians usually use only equivalencies. Definition by abstraction is problematic from the logical point of view[15]. But I want to stress a different point. Even if definition by abstraction were justified logically it would not provide what a mathematician normally looks for. Frege's "direction" (not to be confused with orientation!) is hardly an interesting mathematical notion; this concept might play at most an auxiliary role in geometry and can easily be dispensed with. The idea of a *family* of parallel lines does the same job as Frege's abstract direction but is more convenient and more intuitive. Similarly it is more convenient to think of a natural number as a family of equal "doubles" rather than a unique abstract object. Such abstract numbers would be much like Plato's ideal numbers. Plato certainly had a point in arguing that such things do not belong to mathematics!

Frege would most likely answer that the question of convenience does not matter because his proposal is logically justified and the more traditional mathematical practice and parlance is not. An argument justifying that traditional and convenient practice is given in section **9** below.

## 8. Relative Identity

The Theory of Relative Identity is a logical innovation due to Geach (1972, ch.7) motivated by the same sort of mathematical examples as Frege's definition by abstraction. Like Frege Geach seeks to give a logical sense to mathematical talk "up to" a given equivalence E through replacing E by identity but unlike Frege he purports, in doing so, to avoid the introduction of new abstract objects (which in his view causes unnecessary ontological inflation). The price for the ontological parsimony is Geach's repudiation of Frege's principle of a unique and absolute identity for the objects in the domain over which quantified variables range. According to Geach things can be same in one way while differing in others. For example two printed letters *aa* are *same as a type* but different *as tokens*. In Geach's view this distinction does not commit us to *a*-tokens and *a*-types as entities but presents two different ways of describing the same reality. The unspecified (or "absolute" in Geach's terminology) notion of identity so important for Frege is in Geach's view is incoherent[16].

---

[15] See Scholz&Schweitzer (1935) for a historical survey and further references.

Geach's proposal appears to account better for the way the notion of identity is employed in mathematics since it does not invoke "directions" or other mathematically redundant concepts. It captures particularly well the way the notion of identity is understood in Category theory. According to Baez&Dolan (1998)

In a category, two objects can be "the same in a way" while still being different (p.7)

so in Category theory the notion of identity is relative in exactly Geach's sense. But from the logical point of view the notion of relative identity remains highly controversial. Let $x,y$ be identical in one way but not in another, or in symbols: **Id**($x,y$) & ¬ **Id'**($x,y$). The intended interpretation assumes that $x$ in the left part of the formula and $x$ in the right part have the *same* referent, where this last (italicized) *same* apparently expresses absolute not relative identity. So talk of relative identity arguably smuggles in the usual absolute notion of identity anyway. If so, there seems good reason to take a standard line and reserve the term "identity" for absolute identity.

We see that Plato, Frege and Geach propose three different views of identity in mathematics. Plato notes that the sense of "the same" as applied to mathematical objects and to the *ideas* is different: properly speaking, sameness (identity) applies only to ideas while in mathematics sameness means equality or some other equivalence relation. Although Plato certainly recognizes essential links between mathematical objects and Ideas (recall the "ideal numbers") he keeps the two domains apart. Unlike Plato Frege supposes that identity is a purely logical and domain-independent notion, which mathematicians must rely upon in order to talk about the sameness or difference of mathematical objects, or any other kind at all. Geach's proposal has the opposite aim: to provide a logical justification for the way of thinking about the (relativized) notions of

---

[16] For recent discussion see Deutsch (2002)

sameness and difference which he takes to be usual in mathematical contexts and then extend it to contexts outside mathematics :

Any equivalence relation ... can be used to specify a criterion of relative identity. The procedure is common enough in mathematics: e.g. there is a certain equivalence relation between ordered pairs of integers by virtue of which we may say that *x* and *y* though distinct ordered pairs, are one and the same rational number. The absolute identity theorist regards this procedure as unrigorous but on a relative identity view it is fully rigorous. (Geach 1972, p.249)

## 9. Internal Relations

In (1990) Russell writes:

Mr. Bradley has argued much and hotly against the view that relations are ever purely "external". I am not certain whether I understand what he means by this expression but I think I should be retaining his phraseology if I described my view as the view that all relations are external. (p.143)

In arguing that relations are generally speaking *internal* Bradley (1922) means roughly the following: the relata of a given relation generally speaking cannot be thought of independently of each other and of the relation in question. (So relations, if any, such that their relata *can* be thought of independently are *external*[17].) Bradley makes indeed a stronger claim:

Relations exist only in and through a whole, which cannot in the end be resolved into relations and terms. … The opposite view is maintained (as I understand) by Mr. Russell. But for myself, I am unable to find that Mr. Russell has ever really faced this question (1922, p.127)

As we can see each of the two authors admits that he hardly understands arguments of the other. Since Russell's outright rejection of internal relations they have been under great suspicion amongst Analytic philosophers. Today the neglect of internal relations is not only the consequence of underlying inclinations in systematic metaphysics but also a matter of available

---

[17]Notice that there is no duality between external and internal relations since internal relations are *not supposed* to be defined independently of their relata (which would be an absurdity). See Hylton (1990)

logical means. For the main tradition of (modern) logical systems is developed in keeping with Russell's rejection of internal relations, so one may ask whether or not the standard modern notion of n-placed relation as n-placed predicate can be possibly understood as *internal* relation. Let us see. Consider the standard procedure of interpretation of given relation **R**(*x,y*) in given domain ***D***. Here *x*, *y* are logical variables which take their values among members of ***D*** that is usually thought of as a class of individuals. When *x,y* take values ***a,b*** from ***D*** **R**(***a,b***) takes a certain truth value (usually "true" or "false") just like function *f(x,y)=x+y* takes value 3 when *x* takes value 1 and *y* takes value 2 (and + is interpreted in the usual way). Noticeably the substitution of ***a***, ***b*** for *x,y* is proceeded uniformly for any binary relation and in this sense it doesn't depend on **R**. To put it in other words the substitution is *formal*: one first substitutes ***a,b*** for *x,y* and then looks for the true-value of **R**(***a,b***). So relata ***a,b*** are assumed here independently of **R**. This meshes well with Russell's view according to which all relations are external.

But can **R**(*x,y*) be possibly understood as internal? Consider relation **NEXT**(*m, n*) between natural numbers which says that number *m* is followed by number *n*. Arguably natural numbers cannot be correctly thought of without **NEXT**. This means that this relation is internal. But this claim apparently has nothing to do with the order in which **NEXT** is interpreted: nothing prevents one to pick up numbers 1, 2, substitute them for *m, n* in **NEXT**(*m, n*) and see that **NEXT** (1,2) is true. Thus the logical machinery involved here has no bearing on the metaphysical controversy between the external and internal understanding of relations. So given relation **R**(*x,y*) might be internal as well as external.

However the above argument is not convincing. For it involves interplay between the formal analysis of the concept in question and implicit assumptions made about this concept. As far as we *already know* what are natural numbers then we can claim, of course, that 1 is followed by 2. We can also write down this truth in a more fashionable way as **NEXT**(1, 2). Formal logic is used in this case for description of a ready-made concept. In such a case logic has no bearing on how the concept in question is built, and so it is metaphysically neutral. But when logic is used for *concept-building* like in the case of foundations of mathematics then specific features of logical apparatus get directly involved into emerging concepts. In practice the distinction between the two ways of applying logic can hardly be ever made rigidly: the major application of logic is a logical *reconstruction* of given background (in particular of common mathematical practice) but not an external description of ready-made concepts nor creation of new concepts from nothing. Logical reconstruction is making of new concepts from *something*. I will not elaborate this point here and only notice that foundations of mathematics obviously involve concept-building even if it has a descriptive function too (with respect to the common mathematical practice).

Is it possible to stipulate the relation **NEXT** between natural numbers *without* assuming a full-fledged notion of natural number in advance? A positive answer is given with Hilbert's axiomatic method. One assumes some *class of individuals N* as domain of binary relation (two-placed predicate) **NEXT**($m$, $n$), and stipulates certain formal properties of **NEXT** as axioms. The idea is that a system of axioms of this kind will turn abstract individuals into numbers. Or to put it more accurately, elements of given class *N* will be thought of as natural numbers as far as they verify some properly chosen axioms. Think about Peano's arithmetic.

Is this procedure indeed compatible with the internalist account of relations? The answer is not trivial, in my view. On the one hand, there is obvious reason to think of **NEXT** introduced axiomatically as *internal*: unless **NEXT** (with its formal properties) is taken into account elements of *N* are thought of as abstract "things" but not as numbers. But on the other hand, the stipulation of relata of **NEXT** as individuals (elements of given class) is incompatible with a strong version of internalism about relations according to which these relata cannot be thought of without its relation at all, not even as abstract "things" without properties. So the standard logical apparatus is indeed incapable to represent relations which are internal in this strong sense. Apparently Bradley defends such a strong version of internalism about relations when he says that "a whole … cannot … be resolved into relations and terms". True, this radical position undermines the very notion of relation, so after all Russell's account of relations should be probably preferred. However Bradley's remark points to a real problem which shows that the notion of relation (or at least in its Russell's restricted version) is far less powerful than it seems. Notice that any axiomatic theory based on relations (for example ZF based on membership) assumes its objects (for example sets) as individuals. However we have seen that the identity of basic mathematical objects like points, circles or natural numbers is highly problematic. The blunt stipulation of such things as individuals doesn't resolve the problem but turns it into a new form: given two classes *N* and *N'* (*N'* might be a "copy" of *N*) both satisfying axioms of arithmetic which of the two classes is *the* class of natural numbers? (Benacerraf 1965) I suppose that in order to get a satisfactory solution of the identity problem in mathematics we should give up the idea that mathematical objects always form classes and look for different ways of getting multiple objects into a whole. In the following section I analyze the notion of class and show its limits.

**10. Classes**

Sets of chairs or crowds of people are usually considered as a paradigm cases for our thinking about the notion of *many*. There are different examples though. Think about clouds in the sky or waves at see surface. One can always count persons or chairs or at least in principle so. But one can hardly count clouds and waves. The problem is not that they are too many but that there is no definite criterion for distinguishing one from another. Clouds and waves are certainly many but this kind of many is in general not countable. For a mathematical example think about families of equal numbers or of coincident points: the question of the cardinality of such *multiplicities* (to choose a term for *many* with the broadest meaning) is apparently senseless. In what follows I shall specify the sense of "countable" relevant to this context. I shall term the wanted concept *weak countability* in order to avoid confusion with countability in the usual set-theoretic sense. In his (1903), ch. VI ("Classes"), Russell distinguishes between extensional and intensional "genesis of classes": the former proceeds through the "enumeration of terms" while the former proceeds as follows: one takes a predicate **P**(*x*) and considers class {*x* I **P**(*x*)} consisting of such *x* that **P**(*x*). For example class {1,2,3} can be defined either through the direct enumeration of its elements 1,2,3 (extensional genesis) or as the class of natural numbers smaller than 4 (intensional genesis). According to Russell the extensional genesis of classes through enumeration is possible only when the number of elements (terms) is finite. However Russell claims that this constraint is only "practical" and "psychological" but not logical and theoretical. In particular he says:

… logically, the extensional definition appears to be equally applicable to infinite classes, but practically, if we were to attempt it, Death would cut short our laudable endeavor before it had attained its goal. Logically, therefore, extension and intension seem to be on a par. (1903, p 69)

After claiming the essential equivalence of extensional and intensional viewpoints Russell goes ahead and claims the priority of extension:

A class … is essentially to be interpreted in extension. …. But practically, though not theoretically, this purely extensional method can only be applied to finite classes. … although any symbolic treatment must work largely with class-concepts and intension, classes and extension are logically more fundamental for the principles of Mathematics. (ibid. p.81)

These arguments are not convincing. True, theories often extend domains of possible application of available practical means through relaxing certain constraints. For example, since Ancient times people tend to think about distances between celestial bodies and between pebbles on sand on equal footing. In many cases such theoretical extension works and allows for improvement of existing practical means; in other cases taking practical constrains into theoretical consideration allows for improvement of theories (think about Gauss' work in geodesy which motivated his geometrical discoveries).  However I cannot see how this might help to settle the issue of extension and intension. What Russell says about Death is irrelevant: an immortal god would no better succeed to accomplish the task of finishing enumeration of an infinite series than a mortal human because it cannot be possibly finished. So to the contrary of Russell's opinion, the difficulty of the infinite enumeration is not practical nor psychological but certainly theoretical and logical. Russell refers to the mathematical (Cantorian) notion of infinite set but he misses an essential point of Cantor's invention. In his (1883) Cantor says roughly the following. Count 1,2,3,… This counting never ends - not practically nor theoretically - but we may stipulate a new ideal object $\omega$ as the limit of this process like we stipulate an irrational number as a limit of a series of its rational approximations. Then $\omega$ can be understood as a number of all (finite) natural numbers, and so the talk about the set of all natural numbers becomes reasonable. Cantor proposes here a specific extension of the usual finitary enumeration, and I don't think that the philosophical distinction between the theory and the practice much clarifies this Cantor's proposal. Observe that Cantor's invention has no immediate bearing on the issue of predication, so Russell's idea that a predicate may bring about anything like Cantorian set (remind that Russell doesn't distinguish between sets and classes) is a very strong independent hypothesis. The following development of logic and Set theory imposed well-known constraints upon the use of classes but these commonly accepted constraints, in my view, are not sufficient and somewhat misleading. Set-theoretic "antinomies" including Russell's paradox forced Zermelo (1908) to restrict Russell's "intensive genesis": Zermelo's "Aussonderungsaxiom" allows for "genesis" of set *S* using property **P** only given another set *M* such that **P** is *definite* on *M*, which means that **P** has definite truth-value for every element of *M* (then *S* comprises those elements *x* of *M* for which **P**(*x*) is true). Russell himself completely changes his mind already in (1906a, b) putting forward *No Class Theory* according to which what one needs in logic is only a domain (class?) *U* of individuals but not any further classes constructed out of *U*. In (1908) Russell changes his mind again and puts forward his type theory.

Bernays in his (1958) purports to save Russell's early liberal notion of class through a formal distinction between classes and sets[18]. According to Bernays sets are classes having a specific property of being *individuals*, that is, capable of being elements of other classes. For sets Bernays accepts an improved version of the Aussonderungsaxiom (which he proves as a theorem). However classes in Bernays' view are formed by properties "automatically", so one even doesn't need a quantifier for it and can simply write {*x* I **P**(*x*)} to denote the class of all *x* such that **P**(*x*) is true (this class can be a set or a proper class dependently on predicate **P**).

Moreover Bernays doesn't exclude the possibility that classes can be produced in other ways not mentioned in his theory:

This point of view suggests also to regard the realm of classes not as fixed domain of individuals but as an open universe, and the rules we shall give for class formation need not to be regarded as limiting the possible formations but as fixing a minimum of admitted processes for class formation. (1958, p. 57)

---

[18] The distinction between proper classes and sets has been introduced earlier by von Neumann. See Fraenkel (1958), p.32-33.

Bernays' liberal notion of class remains very popular among mathematicians. People have learnt that the notion of set shouldn't be applied without caution but thanks to Bernays they feel free to talk about classes of anything. This has changed the way of thinking even about elementary mathematical concepts. The idea that the Euclidean plane contains the class of all circles would sound completely weird in the 19th century but today's mathematical students usually don't feel any inconvenience about it. In the eyes of many this freedom of thinking about infinite collections ("Cantorian Paradise") is a very important achievement of mathematics of 20th century.

A standard worry about such extensional representation of mathematical concepts concerns the issue of infinity: why we need such huge collections where we can do well with only few examples? I'm rather sympathetic with the parsimony reason behind this worry but now I want to stress a different point. Another worry concerns the fact that thinking of, say, circles on the Euclidean plane, as forming a class we are obliged to take circles as full-fledged *individuals* with definite identity criteria. But as I have tried to show in the beginning of this paper (section **2**) such criteria are hardly available.

Before I elaborate on this crucial point let me make a methodological remark. When I criticize the class-based representation of mathematical concepts I do *not* assume that there exist the only right way to represent mathematical concepts. I believe that mathematical concepts are exactly what we think about them, and that there is a sense in which the same concepts can be represented differently. The class-based representation is *a* way of thinking about mathematical concepts which proved to be in many ways successful. My critical efforts directed against this approach aim at revealing its hidden assumptions and constraints and at giving place for alternative approaches, which look more promising. When I say "circles are not individuals" I mean that the class-based representation of circles clashes with what people usually think about circles in many standard contexts. I recognize that this clash alone provides no strong argument against the class-based representation: perhaps we should fix the traditional way of thinking about mathematical objects rather than modern formal methods. However in the following sections I shall show that these traditional intuitions support some important contemporary mathematical developments, so in order to promote these developments we need to elaborate on these intuitions rather than rule them out.

**11. Individuals**

Bernays understands the notion of individual in the logical sense as an element of a domain of quantification, that is, an element of some class. The extensionality property of classes (which Russell rightly stresses as indispensable) implies that individuals so understood (elements of classes) must have unproblematic identity criteria. To see this remind how the Axiom of Extensionality is written in ZF:

EXT: $\forall x \forall y (\forall z (z \in x \leftrightarrow z \in y) \rightarrow x=y)$

Informally this axiom says that "sets are wholly determined by their elements". Although the identity of sets is introduced in ZF independently of EXT[19] the intuitive appeal of this axiom certainly depends of the fact that it can be used for "checking identity": given two sets one can check whether or not they are the same through checking their elements. Remark however that in ZF there is no distinction between sets and elements as different types: elements of sets are also sets. So EXT reduces the question about identity of given pair of sets to the question about identity of some other pairs of sets. If given sets $x,y$ are infinite then checking the identity $x=y$ through EXT reduces the problem to checking an infinite number of identities. Prima facie this doesn't look helpful. In fact EXT is helpful for checking identity $x=y$ *only* when the questions about identity of elements of $x, y$ have obvious answers or at least are easier to answer. If identity of elements of $x,y$ is just as problematic as identity of $x,y$ themselves EXT looses all of its appeal.

---

[19] To make EXT into an instance of Keranen's *identity schema* we need to replace the implication by the biconditional: EXT': $\forall x \forall y (\forall z (z \in x \leftrightarrow z \in y) \leftrightarrow x=y)$. EXT' is true in ZF but is not used neither as a definition nor as an axiom for the reason of logical parsimony. In fact ZF allows for another instance of the identity schema obtained from EXT' by the reversal of $\in$: INT : $\forall x \forall y (\forall z (z \in x \leftrightarrow z \in y) \leftrightarrow x=y)$. Taking INT as giving the sense of identity brings about rather unusual way of thinking about sets, which I developed in my (2003b).

In ZF the "identity check" with EXT finishes after a finite number of steps thanks to the Foundation axiom, which forbids chains of the form … $\in x \in y \in z$…. infinite to the left and cycles of the form $x \in y \in … \in x$. So the procedure always ends up with *the* empty set (or should I say with a *copy* of the empty set?). This is hardly what predicts the common intuition about sets based on examples of sets of points and the like[20].

Bernays assumes the extensionality of classes but in order to avoid quantification over classes he modifies EXT into this open formula

EXTCl : $\forall z(z \in x \leftrightarrow z \in y) \leftrightarrow x=y$

which he uses as the definition of identity (equality) of classes (here $x,y$ are classes while $z$ ranges over sets). So the extensionality of classes in Bernays' account becomes also "automatic" and doesn't require a special axiom. Anyway EXTCl provides classes by definite identity conditions just like sets. However according to Bernays certain classes (proper classes) cannot be elements of other classes. *Why not?* Because it is known that making classes elements of other classes in certain cases leads to contradiction. But this is a mere recognition of the fact but not an explanation of the phenomenon. The colloquial explanation according to which proper classes are "too big" or "over-comprehensive" (by Fraenkel's word in his 1958) for being elements of something bigger (because there is nothing bigger?) certainly cannot be viewed as satisfactory.

---

[20] The fact that all sets in ZF are "made of the empty set" (just like Zermelo's or von Neumann's ordinals) is indeed quite counterintuitive. This gives reason for introduction of "atoms" (otherwise called *urelements*) which are elements of sets without their own elements but which are, generally speaking, not all identical to each other. The introduction of atoms restricts the extensionality since EXT implies that sets without elements are identical. A way to introduce atoms is to declare them as things of different type than sets, and to make variable $z$ in EXT to range over sets *and* atoms while leaving $x,y$ to range only over sets. In a set theory with atoms the identity checking with the (restricted) extensionality axiom ends up with the question of identity of atoms, so this question must be unproblematic if we want to respect the basic intuition about extensionality.

Here is my explanation, which implies a substantial revision of Bernays' point of view. I suppose that *multiplicities* like "all sets" cannot be viewed as individuals because *their* elements are not individuals either and hence have no definite identity conditions. Such multiplicities cannot be thought as classes (or as elements of other classes) on the pain of loosing the sense of extensionality. Although we can think about all sets in a way we cannot think of all sets as *individuals*.

I'm not presently prepared to offer a full theory of mathematical individuation but I shall make a hypothesis, which allow me to defend the claim. In the traditional (pre-Cantorian) mathematics the individuation is always finitary and associated with *naming*: one stipulates, for example, points $A$, $B$, $C$,… (some of which might appear to be identical) but not an infinite set of individual points. This doesn't, of course, preclude one of speaking, say, about "any point of given line"; the difference with the modern point of view is that this expression doesn't commit one to an infinite set of points. Cantor's notion of infinite set is based on the assumption that individuals can form not only finite but also infinite collections. In other words he assumes that thinking about *all* points of given line we can still think of these points as individuals like $A$, $B$, $C$. Cantor provides the following justification of this view. He shows that a properly generalized procedure of *counting* (enumeration) of elements of given set works in the infinite case too. (This applies to all infinite sets but not only for sets which are countable in the usual technical sense, see about Bernays' Numeration Theorem below.) This doesn't really prove that elements of infinite sets are individuals in precisely the same sense as elements of finite sets like $\{A, B, C\}$ but this shows that at least one essential feature of finite sets is preserved in the infinite case, namely, the fact that elements of infinite sets may be brought into one-to-one correspondence just like elements of finite sets $\{A, B, C\}$ and $\{D, E, F\}$. This gives indeed a reason to think of elements of infinite sets as individuals by analogy with the finite case.

I shall call multiplicities having a cardinality *weakly countable* and require classes to be weakly countable. Given this additional requirement for classes I shall call elements of given class *individuals*. Thus my hypothesis is that weak countability implies (at least a weak form of) individuation[21].

This hypothesis is in accord with Russell's point that all classes are in a certain sense "denumerable". Unlike Russell Bernays says nothing about enumeration of classes but for sets he proves the *Numeration Theorem* (*ibid.*, p.138) which improves upon Cantor's infinitary "enumeration" in terms of formal rigor and states that every set has a certain cardinality. The theorem doesn't hold for proper classes. Nevertheless Bernays assumes that proper classes consist of well-distinguishable elements, and that the extensionality property holds for proper classes. In my view this assumption is ungrounded. Just like Russell in (1903) he apparently thinks that that a mere predication brings about some sort rudimentary enumeration. I don't think that this is a correct view.

Consider predicate *human* for example. The collective term *humans* unlike *all presently living humans* is not associated with any particular group of people. The expression *all humans* doesn't make much sense unless it is further specified (in any event we cannot count all future generations). Nevertheless we can speak about *humans* as a multiplicity. When we talk about *sets* in mathematics the situation is not different. Multiplicities of all sets or of all singletons don't deserve the name of classes because such multiplicities have no definite cardinalities and hence there is no reason to think of their elements as individuals. (Notice that according to Bernays every set as an element of the class of all sets *is* an individual.) So the problem is not that sets, singletons, etc. are too many but that there is no certain (cardinal) number of such things at all.

Bernays disqualifies Russell's aforementioned definition of cardinal numbers as classes of equivalent sets because he wants them to be sets but not proper classes. He doesn't say explicitly why but apparently he himself has certain misgivings about proper classes and doesn't want to see them playing a crucial role in his set theory. So he identifies cardinal numbers with certain ordinals. This technical solution causes Benacerraf's problem already mentioned: why we should call *cardinal number* one particular set of given cardinality rather than another? This kind of definition of cardinal number differs drastically from Frege's and Russell's earlier proposals discussed in section **4**, so what I mention then as a success of their project should be taken with a pinch of salt.

---

[21] Equating the weak countability with having certain cardinality I take the most liberal attitude possible intended to preserve the whole of Cantorian set theory. More constructively-minded people might prefer to equate the weak countability with the usual countability, or even to insist that infinite enumeration is impossible.

Thus my point is that weak countability required by classes shouldn't be always taken for granted and expected to be found everywhere in mathematics. As the phenomenon of "mathematical doubles" suggests many mathematical objects might be accountable in terms of internal relations (in particular internal equivalences) which don't allow for considering these object as full-fledged independent individuals. Moreover the unique "full-fledged " notion of individual (and hence the unique notion of identity) should be likely given up in favor of various specific structures. (In particular, the kind of individuation implied by the weak countability is obviously weak itself because it allows for arbitrary permutations between "individuals". Given, say, a pair of elements this kind of individuation allows one to say that the elements are indeed two but not to specify which is which.) Such non-standard identity structures can be hardly introduced through standard logico- semantical techniques, which assume the notions of class and logical individual as previously given. Let us look for a different framework.

**12. Relations versus Transformations**

The replacement of the equivalence $x\mathbf{E}y$ by the identity $x=y$ proposed by Frege allows for a stronger interpretation than Frege gave in his account of abstraction. Namely, **E** can be interpreted as a *reversible transformation,* which turns $x$ into $y$ and vice versa, and the identity $=$ as identity *through* this transformation. In the case of congruence the transformation is (Euclidean) *motion*: $y$ is the *same* object $x$ but subject to translation and/or rotation in Euclidean space. Here $x$ and $y$ are said to be the *same* in the same sense in which, for example, an adult yesterday and today is the same person. So we think here of a given triangle in much the way we think of a *substantial continuant* - as an entity capable changing its states and/or positions. Such a "substantialist" interpretation works also for Frege's example of parallel lines. States and positions are *secondary entities* in Aristotle's sense: they exist only in virtue of the (existence of) a substance they are states and positions *of*. To put it in modern terms positions are ontologically *supervenient* upon other entities[22].

The substantialist reinterpretation of mathematical relations may look like an exercise in old-fashioned metaphysics but it appears surprisingly fruitful from the mathematical point of view. For in mathematics the language of transformations is not formally equivalent to that of relations as one might expect but is actually far richer. Given equivalence $x\mathbf{E}y$ there are, generally speaking, *many* distinguishable transformations turning $x$ into $y$ while $x\mathbf{E}y$ only says that one such transformation exists. So here the underlying naive metaphysics matters mathematically.

The difference becomes particularly evident in the case of (global) reversible transformations of a given geometrical space. In the language of relations the existence of such transformations amounts only to the claim that a given space is equivalent to itself. But in fact such transformations contain the most basic information about the corresponding space. This was first recognized by Klein in (1872) when he formulated a new research program in geometry as follows:

Given a manifold[23] and a group of transformations on it one should investigate the structures on the manifold with respect to those properties that respect the transformations of the group. (p.7)

It is not the notion of a "substantial form" surviving through transformations that is the major issue in the new framework for the study of geometrical structure proposed by Klein. Rather there is something of a different sort, which also remain unchanged through the transformations.

---

[22] For a more up-to-date account of the notion of substance and of identity through change see Wiggins (1980)

[23] "Es ist eine Mannigfaltigkeit und in derselben eine Transformationsgruppe gegeben; man soll die der Manningfaltigkeit angehoergen Gebilde hinsichtlich solcher Eigenshaften untersuchen, die durch die Transformationen der Gruppe nicht geaendert werden." Of course Klein here is not using the term "Mannifaltigkeit" in the technical sense which has become standard today.

That something is the structure(s) or *forms* of the transformations themselves. I refer to the fact that reversible geometrical transformations like Euclidean motions form algebraic *groups* under composition. This fact remains completely hidden from view when one uses the language of relations. Thus the traditional metaphysics of substance and form fulfills a mathematical need which the new Frege-Russell metaphysics does not - whatever might be said in favor of the latter against the former for philosophical reasons.

Let me next specify some terminology, which will be useful for what follows. We have considered three different ways of thinking of what is involved in operating with an (arbitrary) equivalence relation xEy.

1) <u>Extension</u>  Consider equivalence classes formed of those things equivalent under the relation **E**

2) <u>Abstraction:</u> Replace the relation $x\mathbf{E}y$ by identity $x=y$, and read $x,y$ anew as standing for a (relational) property common to all and only members of the same equivalence class under **E**.

3) <u>Substantivation</u>. Think of the given relation as a reversible transformation of relata into each other, and read E as identity *through* this transformation.

In the case of Humean relation **H** one may proceed from 1) to 3) through the following steps. Given certain class of classes *x*, *y*, … equivalent by **H**

- think of the one-one correspondences between elements of given classes *x,y* as reversible transformations (isomorphisms) *f,g*... turning elements of *x* into elements of *y* and conversely (reversibility implies that different elements of *x* turn into different elements of *y* and vice versa)[24];

---

[24] Noticeably such a reading is found already in Cantor (1895) where he says: "es verwandelt sich dabai *M* in *N*" to the effect that elements of a source set *M* are replaced one-by-one by elements of an equivalent target set *N*.

- think of *x*, *y* as different states of the same underlying substratum *X*, and think of (auto)morphisms *f*, *g*,... as changes of *X*;
- similarly identify all classes equivalent to *x* and *y* with *X*.

A non-trivial fact, which makes mathematical sense of this metaphysical exercise, is that the automorphisms of *X* form a group called its permutation group or symmetric group. To see better what we gain and what we might lose in switching from relations to transformations consider the following table:

| Extensional reading | Substantional reading |
|---|---|
| we write *x=y* for "class *x* is equivalent (isomorphic) to class *y*" | we write *f:X→X* or simply *f* for an isomorphism from a class *X* to itself (automorphism) |
| = is an equivalence relation; this means that: | Automorphisms of *X* form a group; this means that: |
| = is transitive: *x=y* and *y=z* implies *x=z*. | given automorphisms *f*, *g* there exists a unique automorphism *fg* resulting from the application of *g* after *f*. |
| = is reflexive: every class *x* is isomorphic to itself: *x=x* | there exist an identity automorphism 1 such that $1f=f1=f$ for any *f* |
| = is symmetric: if *x=y* then *y=x* | every atomorphism *f* has an inverse $f^{-1}$ such that $ff^{-1}=f^{-1}f=1$ |

Comments on the table:

<u>line 1</u>: Classes *x*, *y* from the left column are identified in the right column as explained above. Notice that *x=y* is a *proposition* while *f:X→X* is a thing, namely a particular morphism (function). Proposition *x=y* *says that* there exists an isomorphism between *x* and *y*, while *f is* such an isomorphism. Things contain more information[25] about themselves than mere evidences of their existence. Hence not all available information contained in the right column is contained in the left. However things themselves usually do not tell us anything. Notice that if *x=y* does *not* hold we still have a proposition, which tells us something. Such information cannot be provided by means used in the right column. So the translation between the language of relations and the language of morphisms is not completely transparent here in either direction.

<u>line 3</u>: Transitivity of = does not wholly capture the concept of composition of isomorphisms, in particular it doesn't capture the fact that the composition is unique.

<u>line 4</u>: The reflexivity of = does not capture at all well the concept of the identity isomorphism, so the "translation" in this line is superficial and even misleading. For the reflexivity does not imply the distinction between identity isomorphisms and other isomorphisms.

<u>line 5</u>: The reversibility of isomorphisms is not wholly captured by the fact that = is symmetric (notice the reference to the identity in the right column and see the comment about the previous line).

Does the approach outlined above provide any viable alternative to Frege's project of settling the question of identity in mathematics by external logical means? Prima facie it seems that the notion of identity through change (transformation) invoked here remains completely informal and not likely to be helpful in avoiding paradoxes mentioned in section **1**. However I claim we have here a new formal concept of identity as the *unity* of a group of transformations. This group-theoretic notion of identity meshes well with the metaphysical intuition that any changing entity contains a core invariant through changes. Merging equivalence classes *x*, *y*, ... into one "class-substance" *X* indiscriminately, we recover a notion of identity as a particular transformation (and one unique for a given group) which we may speculate is connected with the notion of *repetition*[26]. This group-theoretic identity is obviously relative in Geach's sense (for it depends on the type of transformation involved)[27]. It is not immediately clear whether this group-theoretic identity has anything to do with the logical notion of identity , which was Frege's concern. But at least we get a well-defined identity concept here, and one which makes the metaphysical intuitions behind it precise.

---

[25] I talk about "information" informally here.

[26] that is, repetition of the *same* - but here the notion of repetition is thought of as giving meaning to the notion of "the same" rather than the other way round. Cf. Deleuze (1968).

There are at least three objections, which can be brought against the suggestion that we should take this group-theoretic notion of identity as a serious candidate for the philosophical explication of the notion of identity either inside or outside mathematics.

(i) The logical (and metaphysical) notion of identity should apply to the widest possible domain of entities, so one can say which things in a given domain are the same and which differ. But group-theoretic identity corresponds to ("identifies"?) only one object, namely its group.

(ii) Group-theoretic identity does not allow us to form *propositions* like "*A* is identical to *B*". Generally, the group-theoretic identity 1 like any other element of a given group is a particular mathematical object while identity is a basic logical concept. In the group theory like elsewhere in mathematics the role of logical identity is played by mathematical equality =. For the sake of the argument we can now ignore subtle differences between the logical identity and mathematical equality discussed in the beginning of this paper. Anyway 1 and = have little if anything in common except the common name and some vague metaphysical intuitions behind it. For theoretical reason we need to distinguish the two things sharply and reflect the distinction in the terminology rather then allow ourselves to be led by confusing terms and the vague metaphysics.

(iii) In particular the group-theoretic identity like any mathematical object needs certain identity conditions. These identity conditions matter essentially, for example, when one proves the uniqueness of the identity of a given group. Proposition defines the group-theoretic identity 1 taking the logical identity = for granted. Hence one needs a prior logical notion of identity to cope with the group-theoretic notion of identity, so no way the latter can be a candidate for the normalization or mathematical explication of the more general notion.

(*) there exist 1 such that for any f we have $1f=f1=f$

---

[27] So objects can be identical up to a transformation of one type but different up to a transformation of a different type.

In what follows we will see that these problems can be partly fixed through generalizing the concept of group up to that of *category*.

**13. How to Think Circle?**

In this section I shall summarize what has been already said using an elementary example which later will also help me to introduce categories.

Consider a circle:

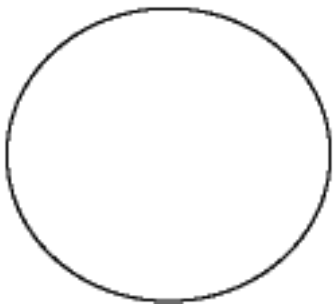

Fig.3

How many things are shown at the picture? We may distinguish at least these two: the general concept of circle *GC* and its instantiation *C*. (Arguably *GC* is shown here *through C*.) Alternatively one may distinguish between circle as a type and *C* as a token of this type (see section **3)**. A different analysis of the same basic situation involves the notions of predicate and individual: think of *GC* as property (one-place predicate) $PC(x)$ and of *C* as individual *IC* such that $PC(IC)$ is true. The "general" character of the concept amounts to the fact that there exist other circles than *C*. This leaves one with a frustrating feeling that there is always something more in the concept than one finds in available instantiations: whatever (finite) number of instantiations one gets some others still remain possible. The Cantorian revolution allows one to overcome this psychological difficulty through bringing *all* individual circles into one infinite class *CC* (see sections **10-11**). Given class *CC* one can observe that all its members are alike (except that they differ in size), and conclude that the worry was basically ungrounded. One could guess this on the basis of few examples, of course. However it would be wrong to conclude from this remark that the multiplying of circles can never bring anything new. Particularly interesting things happen when circles share common geometrical *space* but not just a logical class. Then one can discover that certain circles live inside some others, that some are tangent, etc, etc..

It might be argued that the best way to explain what is circle is to show one as I did it in the beginning of this section. However this would hardly work unless you had seen some circles before. For otherwise you wouldn't know to which features of the shown picture I was trying to attract your attention. A simple way to explain this latter point is to show a learner *many* circles, which differ in size, color, etc., so the learner can grasp common features of these things.

It is also highly recommended to demonstrate to the learner some figures, which look like circles but are not (for example, an oval). After such training the learner normally can correctly identify circles among figures of different shapes.

This pedagogical observation leads to an alternative logical analysis of concepts. Suppose we are given class *F* of figures where circles are mixed up with figures of different shapes. How the circles can be drown out of the mix? If we now are talking in logical but not psychological terms there is the following answer. We need a proper equivalence relation $\cong$ of "alikeness" which would allow us to classify figures from *F* into a number of sub-classes each consisting of "alike" figures, so we would find all circles in one of such sub-classes. Unlike property *PC* equivalence $\cong$ cannot be identified with concept *GC* because $\cong$ helps to distinguish different shapes as well. However for circles $\cong$ does just the same job as *PC*. This shows that one can think circle having nothing like the property of being circle in mind. It may be suggested that property *PC* is construed upon given class *F* and equivalence $\cong$. Remark that $\cong$ can be read as "sameness of shape". Two further alternative suggestions based on this latter remark are available. The first is to introduce shapes as abstract objects in addition to figures, so $\cong$ will serve as the identity relation for shapes (see section **7**). The second suggestion is to diversify *the* identity relation saying that two circles are the *same as shape* but not the *same as figures* (see section **8**). Another way to sort out relationships between a concept and its instantiations is this. Given concept *GC* and class *CC* remark that if circles from *CC* share geometrical space *S* (but not only the class) then any circle can be thought of as a transformed version of any other. So we may replace *CC* by unique circle *SC* which moves around the space and changes its size (but doesn't change shape). These transformations characterizing space *S* and *GC* can be accounted for in terms of Group theory. So we can expect to get an account of *S* and *GC* through a group-theoretic analysis of transformations of *SC* (see section **12**).

## 14. Categorification

Let's look at the method of "moving figures" just described more precisely. The idea that a figure may be moved around within given space (plane) and/or transformed in a way is usual in mathematics. But the idea that there is but one circle on given plane is certainly odd. A mathematician normally needs few such things to play with. Just like in the everyday life in mathematics people deal with many changing and moving things of each kind but not just one. But to the contrary to what one might expect from strict science ambiguities about identity of changing objects seem to be more ubiquitous in mathematics than in the everyday life. One should be a philosopher to wonder if yesterday and today I am the same guy or invent examples like *The Ship of Theseus*. But in mathematics like in the world of Ovidius' *Metamorphoses* (2004) such examples are found everywhere. Given two circles *C, C'* one is always allowed to think of *C'* as a transformed version of *C* making no commitment about preservation of identity through the change similar to that we are usually making in the case of pets and persons. We are so accustomed to these ambiguities through the school education that rarely pay any attention to them. I can quite understand Frege's concern about this issue although I'm not sympathetic with his attempts to improve on it. It might seem that in the end of 19[th] century just like in Plato's times only philosophers could be concerned about the identity issue while mathematicians didn't care. However this is not true. As far as the notion of geometrical transformation becomes basic one cannot always remain unclear about *what* is transformed. Klein suggested this way of avoiding the collapse of circles and other geometrical objects into one: think of transformations of given circle *C* living in space *S* as induced by global transformations of *S* (see section **12**). So ultimately it is the geometrical space but not any particular object in it which undergoes transformations and preserves its identity through the transformations. This assumption makes it easier to treat identities of objects *in* the space in the usual ambiguous way.

However important was the Klein's idea to the time it had a fundamental shortcoming. Isn't it odd to think that in order to move a point one must carry the whole space with it? This idea makes an obstacle for application of new geometrical methods in physics for the simple reason that global transformations of the physical space (or spacetime) have hardly any sense (except in Cosmology)[28]. A purely mathematical aspect of this difficulty is this. Think about a differentiable manifold *M*. Locally *M* is Euclidean, so one has a group of Euclidean motions acting on a neighborhood of each point of *M*. But there is no such group acting on the whole *M* unless *M* is embeddable in Euclidean space of higher dimension. So Klein's approach works smoothly only for Euclidean space which involves no distinction between local and global properties.

Let us see how this problem can be fixed with the notion of category. Instead of choosing between class *CC* of immutable circles and the moving circle *SC* we now allow *all* circles from *CC* to move and transform into each other. This is exactly how things are usually treated in the traditional geometry. More precisely the construction goes as follows. Circles from *CC* are equipped with transformations of two kinds: self-transformations of circles and transformations of circles into other circles. Of course in order to distinguish self-transformations from transformations of the other sort some decision concerning identity of circles must be taken. So there is a space here for a bold stipulation of identity criteria. A natural choice, which suggests itself is this: self-transformations of given circle are rotations around its center.

To complete the definition of *category* of circles **C** there remains only very few things to say.

---

[28] This point has been made by Pierre Cartier during his talk at the seminar on Philosophy of Mathematics in Ecole Normale Superieure in the Spring 2004.

The transformations are composed in the usual way (so the associativity of the composition is assumed) but since the transformations involve different objects one should keep track of what is transformed: transformations *f,g* are composable if and only if the target (domain) of f coincides with source (codomain) of g. Finally with every object *A* we associate a special self-transformation $1_A$ called *identity of A* having the following property: $1_A f = f$ and $g 1_A = g$ for any transformations *f, g* such that compositions $1_A f$ and $g 1_A$ exist (Fig.4).

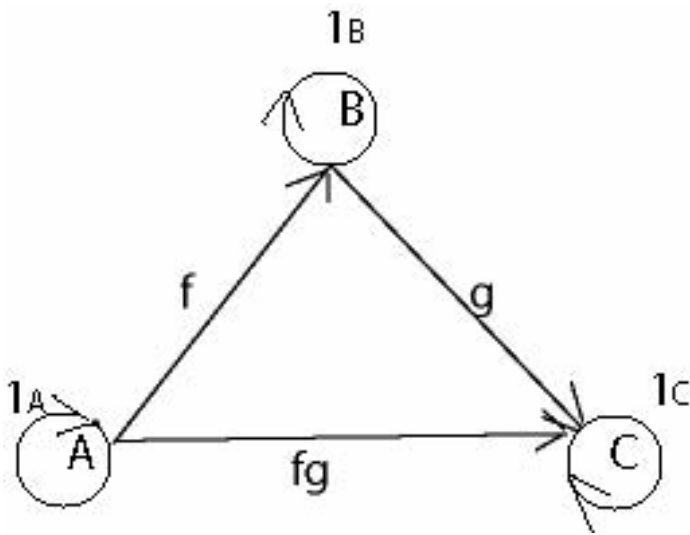

Fig.4

Notice that given circles *A,B* there are, generally speaking many transformations (a class of transformations) from *A* to *B*. In particular there is a class of rotations of given circle *A* transforming it into itself.

In order to get a general definition of category we need only to replace circles by abstract objects and talk about *morphisms* or *arrows* rather then transformations. Thus a category comprises:

- Class of its objects *A, B, C,* …;
- For each ordered pair of objects *A,B* class of morphisms $f: A \to B$, $g: A \to B$ from *A* to *B*; given $f: A \to B$  *A* is called *domain of f* and *B* is called *codomain of f*;
- Composition *fg* of morphisms *f,g* such that the codomain of *f* equals the domain of *g* (see the above diagram); the composition is associative: *(fg)h=f(gh)=fgh*;

- Identity morphism $1_A$ associated with each object $A$ and defined by the condition: for any morphisms $f, g$ $1_A f = f$ and $g1_A = g$ (provided the compositions $1_A f$, $g1_A$ exist).

When in a categorical diagram any arrow $A \rightarrow C$ equals to any other arrow between $A$ and $C$ obtained through composition of arrows shown at this diagram (like at the above diagram) the diagram is said to be *commutative*.

Notice that our category **C** has the following additional property not assumed in the above general definition of category: all its morphisms (transformations) are reversible. The reversibility is a basic property of all usual geometrical transformations like motion or scale transformation in virtue of which such transformations form groups[29]. In the category-theoretic terms just introduced the reversibility of transformation (morphism) $f: A \rightarrow B$ amounts to existence of transformation (morphism) $g: B \rightarrow A$ (called the *reverse* of $f$) such that $fg=1_A$ and $gf=1_B$. In

---

[29] The reversibility of geometrical transformations is a deep intuitive assumption connected to the usual Platonic belief that geometrical objects are eternal (even if not perfectly immutable). Undergoing a reversible change a circle, for example, doesn't take a risk of ceasing to exist and becoming something else. Examples of non-reversible geometrical transformations like shrinking a circle into a point can be easily given but they remain marginal in our thinking about geometrical transformations. To take such examples more seriously would require a significant change of basic geometric intuitions. The reader might object that the importance of circles shrinking into a point and similar things has been recognized in the 19 century the latest through the development of Calculus and particularly in Complex and Vector Analysis. I am agree but I don't consider this as an objection because the way in which such things are dealt with in Calculus (at least from the viewpoint of its official foundations) provides an additional evidence for my thesis. Namely, talking about shrinking of given circle into a point and similar things in Calculus one is obliged to think of an infinite process having a limit which is never achieved. So circle always remains a circle and never becomes a point. Thus the non-reversible transformation of circle into a point, which has nothing mysterious about it and may be provided with many nice properties of typical geometric transformation (like continuity) except reversibility, gets the epistemic status of "naive" and imprecise intuitive picture which may and should be avoided in rigor account.

Category theory this property is taken as the definition of *isomorphism*, so isomorphisms are reversible morphisms. A category like **C** such that all its morphisms are isomorphisms is called *groupoid*. Thinking of objects of a groupoid "up to isomorphism" one gets a group. (So group is a category with only one object such that all its morphisms are isomorphisms.) However such identification causes the lost of information, namely the lost of distinction between morphisms of objects to themselves (automorphisms) and morphisms of objects to other objects. As I have already argued in geometry such a reckless identification of isomorphic objects causes the lost of distinction between local and global properties. Thus groupoids provide an important counter-example against the wide-spread belief according to which in categories all isomorphic objects can be always viewed as identical (see the next section).

The full strength of the notion of category is revealed through the case when morphisms between objects are not all reversible, that is, are not all isomorphisms. A basic example is the category of sets having sets as objects and functions between sets as morphisms. Further examples are obtained through equipping sets with various structures like group structure or topological structure. Then morphisms are required to "preserve" or "respect" the corresponding structure: so in the category of groups morphisms are homomorphisms of groups, and in the category of topological spaces morphisms are continuous transformations[30].

---

[30] Recall the definitions. Homeomorphism between groups $G, G'$ is function $f$ between underlying sets of $G, G'$ such that if $a \otimes b = c$ in $G$ then $f(a) \oplus f(b) = f(c)$ in $G'$, where $\otimes$ and $\oplus$ are group operations in $G$ and $G'$ correspondingly. Continuous transformation $f: T \to T'$ between two topological spaces is a function between underlying sets of $T, T'$ such that any inverse image of an open in $T'$ is open in $T$. Notice that the talk of 'preservation of structure' is at least partly misleading because it too easily makes one think about homomorphisms as if they were isomorphisms. Consider the case of group homomorphism f such that $f(a) = f(b) = f(c) = 1_{G'}$: the "structure" of group $G$ is not preserved in anything like the usual sense of the word but reduced to group unit. The talk of "respect" of structure is less popular but in my view better fits its intended meaning.

Thus the upgrade of the notion of group up to that of category involves two independent steps: (i) introduction of multiple identities (multiple objects) instead of unique identity (unique object), and (ii) allowing for non-reversible morphisms. This upgrade can be shown with the following diagram (Fig.5)[31] :

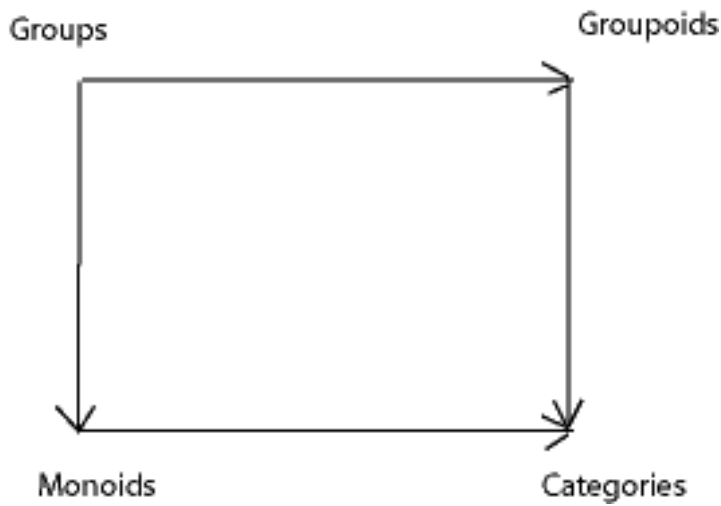

Fig.5

Examples of categories given so far are *concrete* categories. This means that objects of such categories are specified in advance (usually this means that they are construed à la Bourbaki as structured sets), so a category could be seen as a structure over and above given class of its specific objects. However Category theory allows for a different approach: starting with the general notion of category one specifies its algebraic properties to the effect that the structure of morphisms between objects and their compositions determines properties of these objects. The specification of given abstract category amounts to the requirement that certain morphisms exist and certain diagrams commute. In particular a properly specified abstract category "becomes" the category of sets (Lawvere 1964) in the sense analogous to that, in which logical variables in axiomatic systems like ZF "become" sets under its intended interpretation. I cannot provide here

---

[31] In the standard set-theoretical setting monoid is defined as set $M$ provided with a binary operation and unit (identity) 1 ($1 \otimes m = m \otimes 1 = m$ for any $m$ from $M$). Unlike the case of group the existence of inverse elements is not required.

a systematic comparison between the standard axiomatic and the categorical approach to theory-building in mathematics but make only two remarks about this. First, I shall point to what seems to be a difference but in fact is not. At the early stage of Category theory people often opposed categorical foundations of set theory to the standard foundations as "external" approach to "internal". The idea is that while in the case of standard foundations sets are reconstructed though their elements, that is, "from inside", in categorical foundations sets are taken as black boxes "interacting" through morphisms (functions), so what sets *are* is ultimately determined in "sociological" terms of their mutual behavior. This is a right point as far as it concerns basic intuitions about sets but from the formal point of view such intuitions are not essential. The relation of membership $x \in y$ taken as basic in ZF and its likes can be read in both senses – from the left to the right and from the right to the left – and this makes no formal difference (although the intuition behind the extensionality axiom makes the former "internal" reading preferable). At the same time it is not correct that the categorical approach doesn't allow one to "look inside an object": in particular the relation of membership can be perfectly reconstructed by categorical means. In both cases objects of given theory (in particular sets) are first taken as abstract individuals and then "interactions" between the objects tell us "what these objects are". A real difference between the two approaches (and this is my second remark) concerns how exactly these "interactions" are accounted mathematically. In the standard axiomatic approach they are interpreted as relations, and relations in their turn are formalized as predicates (like the two-place predicate $\in$). In the categorical approach the "interactions" are accounted for as morphisms (transformations):

The crystallized philosophical discoveries which still propel our subject include the idea that a category of objects of thought[32] is not specified until one has specified the category of maps which transform these objects into one another and by means of which they can be compared and distinguished. (Lawvere 1991, p.1)

---

[32] Words "category of objects of thought" apparently paraphrase Cantor's famous definition of set given in the beginning of his (1895) as "a collection into a whole of ... objects of our intuition or our thought".

As the example of categorical introduction of sets shows categories represent concepts in a specific way. Thus the mathematical use of the term "category" is not in odd with how this term has been used throughout the history of philosophy.

**15. Identity in Categories: Overview**

The mathematical notion of category just introduced makes paradoxes of identity of mathematical objects discussed in the beginning of this paper more explicit than usual. Consider the category of (all) groups *G*, for example, and take $S_2$ (symmetric group with two elements: unit and involution) as an example of group. Outside *G* one may think about $S_2$ either as a particular object (*the* symmetric group with two elements) or as *kind* of isomorphic objects ("copies") dependently on given context just like one does it with numbers, circles and whatnot. However since category *G* is supposed to comprise *all* groups (whatever this might mean) switching between different senses of "the" cannot any longer remain unnoticed. Similar problems arise in abstract categories. The notion of *terminal object* defined as object having exactly one incoming morphism from each object of given category (including itself) is a typical example. This definition immediately implies that any two terminal objects are isomorphic, and moreover that there is exactly one isomorphism between any two such objects. In any reasonable context (I don't know about exceptions) terminal objects can be identified "up to unique isomorphism", and this is exactly what people do. This identification cannot be hidden by switching to a new context and should be mentioned explicitly. Having no suitable theory of identity in hands category-theorists often justify their liberal use of the equality sign by remarks like this one taken from Fourman (1977). Referring to a formula involving equality the author makes the following reservation:

Strictly speaking the "canonical" isomorphisms…. are necessary (instead of equality – AR)…. Having realized this it is best, in the interests of clarity, to forget them. (p. 1076)

The fact that isomorphic objects are often (albeit not always) regarded as identical in categorical contexts was used by some philosophers as an evidence supporting the claim that Category theory provides "a framework for mathematical structuralism" (see Landry&Marquis 2005 for a recent summary of continuing discussion on this issue in *Philosophia Mathematica*). Mathematical structuralism is, roughly, the view according to which the identity up to isomorphism is the only kind of identity available for mathematical objects. This view squares well with what mathematicians say in informal remarks like the following:

The recursive weakening of the notion of uniqueness, and therefore of the meaning of "the", is fundamental to categorification. (Baez&Dolan 1998, p.24)

or

The basic philosophy is simple: *never mistake equivalence for equality* (ibid., p.46, italic of the authors)

Notice that the "philosophy" suggested by Baez&Dolan here is in accord with my reconstruction of Plato's views given in section **5** and exactly the opposite to Frege's attempts aiming to *strengthening* "the meaning of "the"" in mathematics through (mis)taking various mathematical equivalences for equality (recall the quote from Frege's *Grundlagen* given in section **7**)[33]. I shall not discuss here mathematical structuralism and its relationships with Category theory but remark that mathematically speaking the issue is far from being straightforward. Notice that in the standard categorical setting explained above the identity "up to isomorphism" doesn't apply to all morphisms. To define the notion of terminal object and the very notion of isomorphism (as reversible morphism) one needs to know precisely which morphisms are equal and which are not. So equalities in categories cannot be simply dispensed with and replaced by isomorphisms in any obvious way.

Another part of the same problem concerns isomorphism *of categories*. It has been widely observed that although this notion is easily definable it is quite "useless" (Gelfand&Manin 2003, p.70). Take category *G* of (all) groups for example. Isomorphic copy *G'* of *G* cannot be anything else but the (?) category of groups. But as far as ***G*** is supposed to comprise *all* groups (including all isomorphic groups) the talk of isomorphic copies of ***G*** comprising all these groups once again doesn't make sense (or even is contradictory is "all" is taken seriously). For this reason equivalence of categories is defined as a weaker relation than isomorphism. To give strict definitions we need the notion of *functor*, which is morphism *between categories* respecting the basic categorical structure in the same sense in which homomorphisms of groups respect the basic group structure. Then isomorphism of categories is defined as usual (as reversible functor). To define the *equivalence* between categories we need also the notion of *natural transformation*, which is morphism *between functors* sharing domain and codomain. A natural isomorphism is reversible natural transformation. Now functor $F: A \rightarrow B$ is called *equivalence* if there exist functor $G: B \rightarrow A$ (called *quasi-inverse* of F) such that *FG* is *isomorphic* to the identity of *A* and *GF* is isomorphic to the identity functor of $B^{34}$.

---

[33] Since Frege interpreted this replacement of equivalences by equality as abstraction, this gives an interesting possibility to account for abstraction in terms of *decategorification* introduced in Baez&Dolan (1998) further on.

The equivalence of categories so defined preserves isomorphisms in categories but doesn't preserve identities. This suggests the following view: the "real" sameness of objects in a category is isomorphism but not equality and the "real" sameness of categories is their equivalence but not isomorphism. However we need equality and identity morphisms (in particular, identity functors) in order to define these notions. So a more precise analysis is in order before making any philosophical judgement about identity in categories.

**16. Equality Relation and Identity Morphisms**

First of all notice that the notion of category comprises two very different identity-related elements: the "usual" mathematical equality and identity morphisms. The categorical notion of identity morphism is a generalized version of the group-theoretic identity (the unit of group) discussed in section **12.** (Remind that groups can be viewed as a special case of categories, namely as categories having only one object and such that all their morphisms are isomorphisms.) Let us see how the categorical generalization of the group-theoretic identity amounts to fixing difficulties mentioned in the end of section **12**.

---

[34] For details see Gelfand&Manin 2003, ch. 2 or any introductory text in Category theory.

The first mentioned difficulty is fixed since a category, generally speaking, comprises many objects but not just one. However it remains unclear in which sense if any identity morphisms in a category identify their objects and help to distinguish between them. This is the core problem concerning further difficulties mentioned in section **12**. In the definition of the notion of category objects are introduced as elements of class (that is, as individuals, see section **11**), so they are supposed to be distinguishable from the outset. Then each object is provided with its identity morphism. So from the formal point of view the notion of object as somewhat different from its identity morphism is redundant and may be left out. As far as we are talking about abstract categories nothing is lost indeed if in diagram $A \rightarrow B$ symbols "$A$" and "$B$" are interpreted immediately as identity morphisms. The notion of object is nevertheless useful when we go to examples and want to point to concrete categories. For without the auxiliary notion of object we "loose the substance" and cannot any longer distinguish between, say, sets and functions or between topological spaces and homomorphisms and are obliged to speak about sets and topological spaces as morphisms which sounds odd. The intuitive notion of morphism as transformation also gets lost since we have nothing any longer to be transformed. This makes a metaphysical challenge particularly appealing for a Heraclitean – how to think transformations without transformed substances? – but this doesn't solve our problem which simply reduces to the problem of identification of identity morphisms. In fact we need to identify all morphisms in given category but not only identity morphisms. Like in group theory this is done through the "usual" equality = which is revoked every time when one composes morphisms *f*,*g* and writes *fg=h*. Thus prima facie categorical identity morphisms just like units of groups have nothing to do neither with the logical identity nor with the mathematical equality.

Remind from section **12** that in the group theory the idea to consider identity as transformation and identity as relation on equal footing can be immediately refuted by pointing to the fact that the former is a particular mathematical object while the later is a basic logical concept. In the category-theoretical context this argument doesn't go through straightforwardly. For the argument is based on the assumption that a well-grounded mathematical theory fixes logic first and then develops its specific subject-matter on this basis, so no particular mathematical construction may have any bearing on such a fundamental logical matter as identity. In section **4** we have seen that Frege held this view. I shall call this view *weak logicism* (unlike strong versions of logicism the weak logicism doesn't purport the reduction of mathematics to logic.) The opposite view according to which logic is nothing but a particular application of mathematics and hence cannot be developed independently of and provide foundations for the latter has been held by such otherwise different thinkers as Boole, Brouwer and Poincaré among many others. I cannot find any suitable term for this general view in the history of ideas (which wouldn't refer to something more specific) and I shall call it (logical) *mathematism.* The exact meaning of the term crucially depends on how one distinguishes between logic and mathematics but in the present discussion I shall not purport to be precise at this point.

Logical mathematism moved the modernization of logic in $19^{th}$-early $20^{th}$ century in its early stages (Boole). However later Frege, Russell, and their followers turned the mainstream development in the field from the mathematisation of logic toward a logical regimentation of mathematics with new logical means. Mathematicians had no concurrent large-scale project to offer at the time, so they kept stressing the importance of mathematical intuition[35] and pointing to the fact that the modern symbolic logic involves quite a bit of mathematics. However Category theory gave mathematists a new opportunity. One may still consider category of sets *S* as a semantic for classical logic CL but one can also reverse the perspective taking *S* as a basic construction and looking at CL as a specific structure associated with *S*. Choosing the latter line one may study "logical properties" of mathematically important categories together with their algebraic and geometrical properties and reasonably distinguish between these kinds of properties in categorical terms[36]. The usual distinction between logical syntax and semantics blurred by this approach can be reproduced within a pure categorical setting where some categories play the role of theories and some other categorical constructions play the role of their models. This allows for a mathematist to argue that logical structures are structures of a special kind which don't enjoy any particular conceptual autonomy and/or epistemic priority with respect to the core of (categorical) mathematics. Foundational problems and related philosophical questions remain open but now mathematists exchange their defense for offence putting forward their own agenda in foundations and organization of mathematics and logic.

---

[35] That is why in the $20^{th}$ century logical *mathematism* is closely associated with mathematical *intuitionism.* I deliberately distinguish between the two.

[36] According to one definition a *logical category* is a Cartesian category with (a) images stable under pullbacks and (b) finite sups of subobjects of a given object which are stable under pullbacks. A *logical morphism* between logical categories is a functor, which preserves finite co-limits, images and finite sups. (Kock&Reyes 1977, p. 295-296) The two definitions taken

I don't need refer to further details of this ongoing project to make my present point: the fact that identity morphisms in a category are elements of mathematical construction doesn't imply that these things cannot perform a logical function. For we know that elements of certain categorical constructions do play a logical role, for example the role of quantifiers and truth-values (Kock&Reyes, Fourman 1977). So we cannot rule out on a general ground the possibility that logical notion of identity could be thought of as a categorical morphism. Let's see what is actually going on.

With a suitable category *T* (noticeably with a topos) one may associate logical calculus *L* called *internal language* of *T* to the effect that each formula provable in *L* corresponds to certain commutative diagram in *T* (soundness). *T* is a semantic for *L* in the usual sense but *T* also represents such features of *L* that the usual (Tarskian) semantic doesn't, for example, the truth-values. This gives reason for the aforementioned reversal of the usual point of view on semantic and syntax and explains the term "internal"[37]. *L* brings with itself identity predicate ≡ while the construction of *T* comprises "external" equality = from the outset. The "adjustment" of *L* to *T* makes ≡ and = basically interchangeable. However this doesn't mean that = (or ≡) gets

---

together define what is a logical property "from below" and "from above": the first tells us which properties are necessary for a category in order to be counted as logical while the second tells us that *only* properties preserved by logical morphisms are logical.

[37] In fact *T* and *L* can be adjusted to each other in such a way that the distinction between the two becomes redundant; in this case it is appropriate to call category *T* *syntactical*. For formal introduction of internal language see Fourman (1977) or McLarty (1992).

"internalized" in the same sense in which people speak about internalization of truth-values and logical connectives: the internalization of logic in a category amounts to representation (if not replacement) of the usual logical syntax by categorical constructions while = is *not* a categorical construction but the "usual" god-given mathematical equality! Identity morphisms of **T** are not used for representing ≡. Thus the standard "internalization of logic in a topos" with an internal language has no bearing on the identity issue. We see that although the idea to account for identity in categorical terms cannot be ruled out on a general ground the standard device of "internal logic" doesn't allow us to realize this idea. Let's look for different possibilities.

## 17. Fibred Categories

The following discussion is based on Bénabou (1985). The idea is the following. Recall that categories have been introduced in section **14** as *classes* of a certain kind. In fact in order to work categories one needs some further classes like the class of all morphisms with fixed domain A, for example. Which properties of classes are used in the "naive" Category theory? Let category **C** be our "object of study" and category **B** be our "optical instrument" for studying **C**. **B** can be thought of as category **S** of sets however we can also consider different possibilities, in particular abstractly defined toposes. Following Bénabou I shall call objects of **B** *sets* (remembering that they could be somewhat different than usual sets) and call classes of morphisms or objects of C *families*. (In what follows *families* will reappear as multiplicities of a different sort than classes.) Now given set *I* (an object of **B**) we may master category **C(I)** called *fiber over I* whose objects and morphisms are families of objects and morphisms of **C** *indexed* by elements of *I*, that is, families of the form $X = (X_i)$ and $f = (f_i : X_i \to Y_i)$ where $i \in I$. Bénabou remarks that speaking about categories naively we assume more than this, so we cannot just fix some sufficiently large set I and use it for indexing every time when this is needed. Namely, we also assume the possibility of *re-indexing*: given families $X=(X_i)$, $Y=(Y_j)$ in **C** where $i \in I, j \in J$ and morphism $u: J \to I$ in **B** we assume that family of objects $X_{u(j)}$ and family of morphisms $f = (f_j: Y_j \to X_{u(j)})$ is uniquely defined and "behaves properly". This allows us to extend **C(I)** through introducing new category **Fam(C)** of *families of C* where objects are families of objects of **C** indexed by different sets and morphisms are pairs of the form $(u,f)$ where $u$ and $f$ are as just described. Morphisms of the form $f=(f_i : X_i \to Y_i)$ we identify with $(id_I , f)$ where $id_I$ is identity morphism of $I$ in $B$. The composition of morphisms in *Fam(C)* is defined in the obvious way. We equip the construction with projection functor $p_C$ which sends every family of objects of $C$ to the set by which this family is indexed and every morphism $(u,f)$ between families to morphism u between sets: $p_C: (X_i) \to I$; $(u,f) \to u$.

Now suppose that we know what equality is in *Fam(C)* and in $B$ but not in $C$. This implies that we cannot think of families (of objects and morphisms of $C$) extensionally as usual. In particular given family $X=(X_i)$ where $i \in I$ and morphism $u: J \to I$ in $B$ we cannot define another family $Y=(Y_j)$ by saying that $Y_j = X_{u(j)}$ because the latter equality doesn't make sense for us. Nevertheless we can achieve the same effect through requiring certain properties of *Fam(C)* and $p_C$ as we shall now see. What we need for it is to characterize morphism $\varphi_{(u,X)} = (u, (id_{Y_j}))$ in *Fam(C)* without using equality in $C$ (($id_{Y_j}$) is the family of identity morphisms of objects $Y_j = (X_{u(j)})$ in $C$). Given $u: J \to I$ and $X=(X_i)$ $\varphi$ is characterized up to unique isomorphism by the following property (i): for any morphism $\psi =(v,g)$ with codomain $X$ in *Fam(C)* and any $v'$ such that $v=uv'$ in $B$ there exist in *Fam(C)* a unique $\psi'$ such that $\psi = \varphi \psi'$ and $p_C(\psi')=v'$. In addition morphisms of the form $\varphi=(u, (id_{Y_j}))$ satisfy the *functoriality* conditions (ii): $\varphi_{(u,X)}= id_X$ (identity morphism of family $X$), and $\varphi_{(uv,X)} = \varphi_{(u,X)} \varphi_{(v,Y)}$ for each $v:K \to J$. Now we use these properties as definition of abstract functor $p:F \to B$ called *fibration over B* (or *fibred category over B*) in the case when only the property (i) is taken into account, and called *split fibration over B* in the case when in addition for each pair $(X, u: J \to p(X))$ one makes a particular choice of $\varphi_{(u,X)}$ (called *splitting*) satisfying functoriality conditions (ii). Thus equality in a category can be defined as splitting of fibration over an appropriate base. Noticeably given a fibration its splitting might not exist or be not unique. I refer the reader for further details to Bénabou (1985).

Bénabou's theory of equality in categories allows for regarding objects and morphisms of given category as families rather than bold individuals (see section **7** about the idea of family); these families can be occasionally split into elements through a (split) fibration in different ways dependently of the choice of base. Such splitting is reverse operation with respect to the informal identification of isomorphic objects and morphisms mentioned in the section **15**, and unlike the latter it is performed more rigorously and more "categorically". This reversal is remarkable: it shows that given the definition of equality through split fibration families are no longer thought of as extensional multiplicities, that is, as classes. Recall however that given fibration $p:F \to B$ categories **F**,**B** are construed in the usual way and in particular assume the "usual" equality of morphisms and objects, so the internalized equality relates only to hypothetical category **C** such that *F=Fam(C)*. This situation is quite analogous to that in the Model theory when a formal theory is interpreted in a semantical structure construed independently of this theory either informally or with the help of a meta-theory. As Bénabou stresses in the end of his paper such "meta-equality" is indispensable "unless you do something different from Category theory". In the following section I shall try to do "something different" but now let us consider another approach to internalization of identity in categories.

## 18. Higher Categories

Given abstract category **C** consider class **Hom(A,B)** of its morphisms *f*, *g*, … of the form $A \to B$. Then turn **Hom(A,B)** into a new category formally introducing morphisms of the form $\alpha: f \to g$ (that is, morphisms between morphisms of **C**). Do this for all pairs of objects of **C**. Observe that morphisms $\alpha$ can be composed in two different compatible ways shown at the below diagram

(Fig.6):

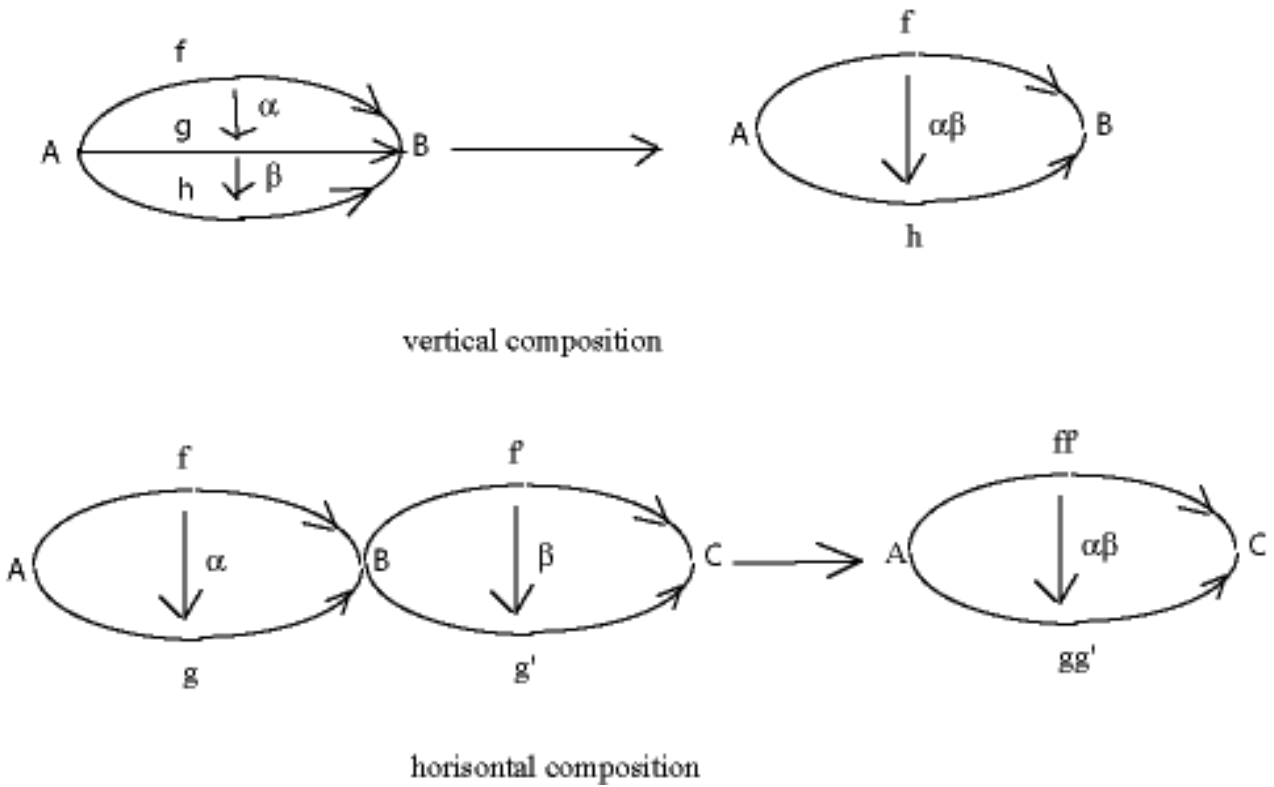

vertical composition

horisontal composition

Fig.6

Requiring now natural equational condition (which I shall not list here) we obtain a *2-category*. It comprises objects $A, B..$, morphisms $f, g,..$ between objects (the same as in **C**) called in this context *1-morphisms*, and "morphisms between morphisms" $\alpha, \beta, ..$ called *2-morphisms.* An example of 2-category which has been around from the very beginning of Category theory (that is, some 20 years earlier than the abstract notion of 2-category has been introduced in Ehresmann (1965)) is 2-category ***2-Cat*** having (some or all) categories as objects, functors between these categories as 1-morphisms and natural transformations between the functors as 2-morphisms. Let me now explain what 2-categories have to do with the internalization of identity (equality). Remark that in a 2-category we have not only the usual composition of 1-morphisms $(f:A \to B)(g:B \to C) = fg: A \to C$ but also functor ***F***: ***Hom(A,B)xHom(B,C)*** $\to$ ***Hom(A,C)*** (provided that in category ***Hom*<sub>C</sub>** having Hom-categories of ***C*** as objects Cartesian product **x** is available [38]).
On 2-morphisms this functor acts as their *horizontal* composition (while in Hom-categories 2-morphisms are composed *vertically*). If functors of the form ***F*** preserve identities in Hom-categories (2-identities) then equalities in ***C*** may be omitted without any lost. This means that we don't even need to define ***C*** as a category but may think of it as a class of "objects" and "morphisms" between these objects, and then define composition of these morphisms "from above" through functors like ***F***. In this case one may speak indeed about "replacement of relations by morphisms": 2-identites from Hom-categories make in ***C*** the job of equalities. The situation here is quite analogous to one we've seen in fibred categories: at the top "meta-" level of construction (namely in Hom-categories and in the category ***Hom*<sub>C</sub>** of the Hom-categories) one uses the "god-given" equality but at the bottom level equalities are got rid of.

An apparent difference between the two approaches is this: in higher categories the notion of *class* is used at all levels including the lowest one while in fibred categories this notion is used only at the "meta-" level and at the lower "internalized" level classes are replaced by non-extensional *families*. But is the assumption that objects of 2-category form a class necessary? Prima facie it is the case. For in order to compose 2-morphisms $\alpha$, $\beta$ in a Hom-category (that is, vertically) we need to check that domain of $\alpha$ equals codomain of $\beta$. So 1-morphisms (objects of Hom-categories) should form a class of well-distinguishable elements (provided with a notion of equality allowing for distinguishing them). Notice that if we take this view then the internalization of equality in *C* just described will be only partial: it will apply to equalities of the form *fg=h* but not to equalities of the form *f=f*. However it is easy to get around this point through identification of 1-morphisms with 2-identities, so all needed equational conditions could be written in terms of 2-morphisms. Then one may claim that in *C* "there is no equality", and hence its elements don't need to form a class.

---

[38] Given objects *A*,*B* in a category their product *A***x***B* is defined up to unique isomorphism by the following (universal) property: given any object *X* with morphisms *X*→*A*, *X*→*B* there exist unique morphism *X*→*A***x***B* such that the following diagram commutes:

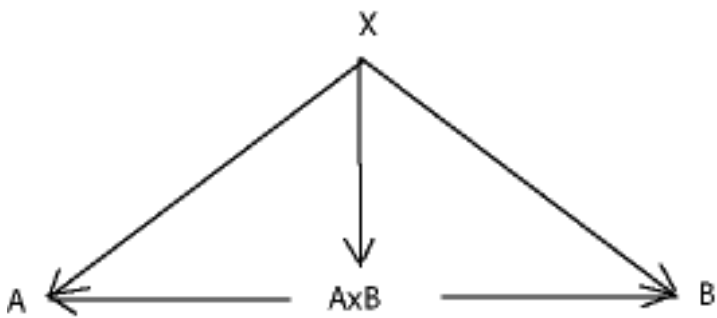

It should be noted that the common interest of people working in higher Category theory is not internalization of equality as such but *weakening* of equality, that is, finding a rigorous way of "replacing» equalities with certain isomorphisms. This approach is quite natural in the given context since the requirement that functors of the form F preserve all 2-identities is unreasonably strong (see section **15**). Since we are no longer obliged to think of *C* as a category in the usual sense we get a room for playing. Instead of imposing on Hom-categories and on **Hom**$_C$ equational conditions implied by the assumption that *C* is a category we can use weakened conditions which don't imply that 2-isomorphisms replacing equalities in *C* are identities. Such weak 2-categories have been first introduced in Bénabou (1967) under the name of *bicategories*. In bicategories (i) the usual associativity law *(fg)h=f(gh)* in *C* is replaced by the requirement of existence of *associativity* (2-)*isomorphism a*: *(fg)h* → *f(gh)* (eventually called *associator* by other authors), and (ii) the usual axiom of identity $1_A f = f$ and $f1_B = f$ for *f:A* → *B* is replaced by the requirement of existence of *unit* (2-) *isomorphisms* *l*:$1_A f$ → *f* and *r*: $f1_B$ → *f*. These isomorphisms are subjects of equational conditions called *coherence laws*, which I shall not list here[39].

The notion of 2-category allows for a straightforward geometrical analogy: think of objects as points, of 1-morphisms as oriented lines, and of 2-morphisms as oriented surfaces bounded by the lines. Whether we use this analogy (which is more profound than it might appear at the first glance) or not the notion of 2-category calls for the inductive generalization to the notion of *n*-category for arbitrary n and further to $\omega$-category (leaving bigger ordinals apart). The *strict* (meaning *non-weakened*) versions of the notions of *n*- and $\omega$-category look unproblematic: the *enrichment* of given category *C* bringing the notion of 2-morphism explained in the beginning of this section can be easily reformulated as an inductive step bringing the notion of *k*-morphism provided with *k* different kinds of composition (Leinster 2001, p.8). However the notions of *weak n-* and $\omega$-categories which are more interesting (both purely mathematically and for applications) are much less obvious. There are many alternative definitions of weak *n-* and $\omega$-categories around; ten of them are presented in (Leinster 2001). A somewhat different approach based on introduction of a new kind of morphism rather than relaxing axioms is presented in Kock (2005). A specific obstacle for putting these things into order, which is not irrelevant to the issue discussed in this paper, is this: it is not clear which notion of equivalence one should apply to answer the question whether two given definitions of weak n-category are equivalent or not. There are obvious reasons to think that the suitable notion of equivalence should be n-categorical itself but this makes the reasoning circular. It seems that new conceptual inventions and not only technical developments are still wanted in this field.

---

[39] See Bénabou (1967) or Leinster (2001), p.9

The construction of *strict n-* ($\omega$) category is transparent because it is built in the strict upward order. However as we have seen in the case of weak n-categories the opposite downward order of construction gets involved. Remind from sections **14** and **17** that the idea of the downward construction is basic for Category theory but not specific for its higher-dimensional branch[40]. At the same time the core upward inductive construction of *n*-category (strict or weak) starting with a class of objects, morphisms between the objects, morphisms between these morphisms, etc., remains present in all alternative definitions of this notion. Although the co-presence of these two opposite orders is ubiquitous in mathematics (think about Euclidean triangles for example[41]) it seems that in the case of *n*-categories the two orders fail to match each other. As far as we start to build a *n*-category from *classes* of objects and morphisms we already assume a lot about their identity, so the sense of further "weakening of identity" through higher-dimensional conditions becomes unclear. For this and other reasons it seems interesting to figure out a categorical (or categorical-like) construction which would avoid the classical upward geometrical concept-building leaving more room for the downward approach. In the following section I shall make an informal suggestion pointing to such a possibility.

---

[40] The talk about dimensions refers to the geometrical aspect of the notion of *n*-category modestly mentioned above in the main text as an "analogy".

[41] I mean the following. To *construct* a triangle one proceeds in the bottom-up order: one takes three points and connects them by straight lines. But to *define* the notion of triangle, say, as a 3-angled polygon one proceeds in the opposite sense assuming the notion of polygon first.

From an epistemological viewpoint the case of *n*-categories with *n*>2 is interesting because it no longer allows for thinking of different levels of the construction along the distinction between the "object level" and the "meta-level". For when n is big the reiteration of "meta-" becomes useless and with $n=\omega$ it becomes senseless. So we cannot take refuge at the "meta-level» but must revise our understanding of identity from the outset. I cannot purport any technical discussion on higher categories in this paper but give the following simple example.

Remind the "partial categorification» of an isomorphism class of sets which we have achieved in section **12**: given such a class we considered all isomorphisms between its member-sets, then identified the sets and some (but not all) isomorphisms and got a symmetric group. Let it be finite symmetric group $S_N$ for simplicity. Now we can see that a more profound alternative to Fregean approach requires turning all sets into a category. Let us however for the sake of example improve on $S_N$ in a different way. Namely, let us categorify it further taking into account isomorphisms of $S_N$, that is, group *Aut($S_N$)* of *automorphisms* of $S_N$. Here we can remark something interesting. Except the trivial cases *N*=1 and *N*=2 when there exist only the identity authomorphisms, and except the "pathological" case *N*=6 we have *Aut($S_N$)* = $S_N$ [42]. (This latter equality sign can be read as the isomorphism relation. Considering isomorphisms in question explicitly we get $S_N$ back!) So taking into consideration automorphisms of higher order (*Aut(Aut($S_N$))* = *Aut$^2$($S_N$)* and so on) brings nothing new: we have *Aut$^n$($S_N$)* = $S_N$ for all *n* and all $N \neq$ 1,2,6. Remark that $S_N$ equipped with *Aut$^k$($S_N$)*, *k*=1,2,…*n* is a very simple albeit not completely trivial example of strict *n*-category. The property of symmetric groups just mentioned is a case of what Baez&Dolan (1998) call *stabilization* in *n*-categories (p.13).

Now we can fix some *n*, assume the "usual" equality only in *Aut$^n$($S_N$)* and for *k*<*n* write the group operation as *ab* → *c* (instead of *ab=c*) claiming that this operation is determined up to isomorphism by the fact that *Aut$^k$($S_N$)* (in particular, *Aut$^0$($S_N$)*= $S_N$) is isomorphic to *Aut$^n$($S_N$)*. The elementary character of this construction makes its downward determination just as easy as the upward one. However it also makes the whole idea of the reiteration of "levels" plainly redundant: given $Aut(S_N) = S_N$ our "n-group" is just $S_N$ "up to itself"! This suggests a different move: instead of describing the group operation through equational conditions reverse the optic and reconstruct a notion of equality (identity) on the basis of the operation. What we can expect to get in this way is "identity up to $S_N$" rather than a universal identity concept suitable for all mathematical and logical needs. Remark however that identity "up to symmetric group of isomorphisms" applies to such a fundamental mathematical object as a *class*. Thus this particular kind of identity has a very general significance in mathematics and logic[43].

---

[42] Kurosh (1960), p.92

### 19. Platonic, Democritean, and Heraclitean Mathematics

It might be argued that before a new account of identity is well established in mathematics it is premature to start any philosophical discussion about it. I don't think so. I dare to think that not only philosophers can find a lot of interesting stuff relevant to their subject in the contemporary mathematics but that mathematicians too can be motivated by new philosophy in their work, moreover if it concerns such a traditional philosophical issue as identity. This is how things worked for centuries (including the heroic time of debates on foundations of mathematics at the edge of 19th and 20th centuries), and I cannot see any reason why in 21st century they should be different. Although philosophical motivations could and in many contexts certainly *should* be swept out of ready-made mathematical theories philosophical reasons often play an important role in bringing new mathematical theories about. That is why the fact that the issue of identity in categorical mathematics is not yet well settled in mathematics gives me reason to discuss it right now rather than to the opposite.

---

[43] I realize the risk of drawing any general conclusion concerning *n*-categories on the basis of the example of "symmetric *n*-group" and suggest the reader to consider this example on its own rights.

In my (2003a) I argued that the classical mathematical tradition dating back to Euclid, which we have all inherited, is profoundly Platonic[44]. However some traces of different mathematical conceptualizations linked to different philosophical stands known in Antiquity are also available. In particular Aristotle's repudiation of thinking about a line as a set of points (*Phys.* 231a*)* shows that the idea was definitely around. Whether or not we can speak about the presence of this atomistic line (which I shall label *Democritean*) throughout the history of mathematics, it is clear that until 20$^{th}$ century it remained marginal. The Democritean mathematics flourished only in 20$^{th}$ century with the rise of Set theory and the triumph of atomistic hypothesis in physics[45]. The modernized atomism of 20$^{th}$ century unlike its Ancient predecessor is "without atoms": our physical atoms and elements of sets are not indivisible but composed of their own elements. What counts here as "atomism" is the very idea that physical and mathematical entities are aggregates of other entities, which are in some sense independent from the former. However as a matter of fact the traditional version of atomism assuming indivisible atoms (points) has been widely revived in mathematics of 20$^{th}$ century as well, so it is common today to think after Cantor, Bourbaki and their unknown Ancient predecessors about lines or topological spaces as sets of *points* provided with an additional structure[46].

---

[44] I'm inclined to explain this by the fact that Plato and his followers used their contemporary mathematical tradition for shaping their philosophy rather than by the assumption that the major part of Ancient mathematicians shared Plato's philosophical views. In some important documented cases the influence was mutual (Fowler 1987.)

[45] In the beginning of his famous lectures (1963) Feynman describes the atomic hypothesis as the most informative claim about the world ever made (p. 1-2).

[46] Today one may always claim in Hilbert's vein that talking about points he doesn't make any commitment about their indivisibility, so points could be thought of as beer mugs or, say, topological spaces (see Hilbert 1899). However unless a non-standard interpretation of points is explicitly mentioned and proves to be non-trivial such getting around doesn't change anything. We can indeed think of point of given topological space as another topological space however the

The assumption about *point* as basic geometrical concept is common for Platonic and Democritean geometry. What differs the two is the way in which points generate the rest of the geometrical universe. In the Platonic case this is the process called in Proclus (1873) $\pi\rho oo\delta o\varsigma$ ("development") which mathematically can be best thought of in *phoronomic* terms (as suggested by Proclus himself): think of line as trajectory of moving point, of surface as trajectory of moving line, and of 3D manifold as trajectory of moving surface. In the Democritean case geometrical objects are constructed with points like with bricks. There are many reasons – purely mathematical and not - why this common assumption is no longer satisfactory[47]. So people are looking for different possibilities.

As far as Category theory is concerned there are two major philosophical motivations behind it. The first is Platonic intuition of mathematical *form* (MacLane 1986). Consider, for example, the standard Bourbakist definition of group as "set equipped with the group structure". This definition makes one think about the background set as a kind of "stuff" and about the structure as a "form" given to this piece of stuff. So the old Plato's metaphor of pottery is still working well here. Then one notices after Plato that what counts here is the form rather than the stuff. Category theory supports this Platonic intuition in the following way. Make groups into a category, forget about the set-theoretic definition of its objects and describe instead the category through equational conditions (commutative diagrams). It looks like we got pure forms without any kind of underlying stuff. When it comes to details the situation appears to be more complicated, and in fact the pottery pattern reappears within Category theory itself (recall, for example, that the categorical notion of fibration introduced in section **17** assumes the "background" category of sets or a background topos).

The other major philosophical intuition behind Category theory I shall call *Heraclitean*.

---

two spaces would have no geometrical relation to each other, and so wouldn't be parts of the same geometrical construction.

[47] See Smolin (2000) for reasons concerning theoretical physics (the assumption that spacetime "consists of points" brought into General relativity through its mathematical apparatus is an obstacle for a theory of Quantum Gravity). For pure mathematical reasons see Cartier (2001).

Many people involved into the development of Category theory perceive the set-theoretic mathematics as somewhat "static" and are motivated by the idea of bringing more "dynamics" in it. I pay more attention to this latter trend in this paper not because of my metaphysical and epistemological preferences but because I believe that the Heraclitean mathematics may open genuinely new possibilities not yet explored, and in particular allow for geometrical frameworks, which are not point-based[48].

Let me now propose a Heraclitean revision of Platonic and Democritean (set-theoretic) ways of thinking about categories. The fact that that Category theory in its present form leaves a space for such a revision seems me obvious. Even if objects of a category are not thought of as structured sets the notion of *class* of objects and morphisms continues to play a principle role in the usual notion of category. In n-categories we have also seen a revival of Euclidean pattern of building mathematical constructions when one starts with points and then continues to lines, surfaces, etc. I shall purport a reconstruction of the notion of category from a "Heraclitean" background without using classes.

For this reconstruction I shall use a mathematical hint borrowed from Euclid's *Elements*. Euclid gives two different definitions of point (equally useless in the Hilbertean view): he defines *point* first, as partless geometrical object ("geometrical atom") (Def. 1.1), and second, as boundary of line (Def. 1.3). He gives similar double definitions of *line* (1.2 and 1.6) and *surface* (1.5 and 11.2)[49]. Aristotle in (*Topic* 141b) qualifies definitions 1.1, 1.2 and 1.5 as "natural" (given "by nature") and definitions 1.3, 1.6 and 11.2 as given "for us". The idea is that while the "natural order of things" requires to start with points and continue to lines, surfaces and solids, for explanatory purposes this order may be reversed, so one may start with solids, then introduce surfaces as boundaries of solids (11.2), then lines as boundaries of surfaces (Def. 1.6), and finally points as boundaries of lines (Def.1.3). This latter order of definitions, which Aristotle qualifies as merely explanatory, is interesting for few reasons. First, if taken seriously, it challenges the common Platonic and Democritean principle of "starting point": this is the reason why Aristotle, who remains a Platonic as far as the issue is concerned, gives to this order of definitions an inferior epistemological status. Second, the "reversed" definitions are perfectly interpretable in terms of today's notion of topological dimension, and so from today's point of view they look more sound than the "natural" definitions based on metric considerations. Third, remark that Aristotle like anybody else until the end of 19$^{th}$ century assumes that all thinkable geometrical objects can be classified into these four types: points, lines, surfaces (including all plane figures) and solids. In the 20$^{th}$ century we have learnt to proceed recursively in the "natural" order and make geometry in higher dimensions. However we haven't learnt yet to proceed recursively in the reverse topological order and still assume that cutting given geometrical object, then cutting the cut, etc., we always end up with an indivisible point after a finite number of steps. Why not to try to go further?

---

[48] The reader may notice that my attitude towards metaphysical matters taken here is pragmatic: I assume that any metaphysics is good if it helps to bring new mathematics about. I realize that this attitude is questionable. Anyway in this paper I leave the question which metaphysic (if any) is "correct" aside.

[49] "Def. *m.n*" refers to *n*-th definition of *m*-th Book of *Elements*

Let's now do some simple mathematics taking the Heraclitean notion of Flow (rather than Point, Solid or Form) as basic. In order to visualize it you might think about flame, boiling water and similar things. People habituated to computer video more than to printed books shouldn't have problems with this. Remark that Platonic and Democritean basic metaphysical pictures can be viewed as two different (and not incompatible) ways of taming the Flow. The Platonic way is stipulation of Forms remaining invariant through the Flow (including Point as a basic Form). The Democritean way is modeling the Flow with a flow of atoms. In both these cases taming of the Flow is conceived not as a mathematical operation but as a metaphysical condition making mathematics possible. I accept the traditional assumption according to which in order to make mathematics with the Flow it must be tamed in one way or another. However I shall tame it through a mathematical rather than metaphysical procedure.

The next primitive notion we need is that of *cut*. Think of cut as a local operation applied to the Flow. Assume that cut $C$ brakes given flow $F$ into four "parts": *incoming* flow $F_C^+$, *outgoing* flow $F_C^-$, *boundary* flow $F_C$ which is both incoming and outgoing, and finally *neutral* (with respect to

C) flow $F_C'$ which is neither. The idea is that while $F$ has no global direction such a direction can be specified *locally* through making local "section" $C$ which allows for distinguishing what "pours in $C$" ($F_C^+$), what "pours out of $C$" ($F_C^-$), and what "flows through $C$" ($F_C$). $F_C$ describes the "internal dynamics" on $C$. Since cut $C$ is local but not global I denoted by $F_C'$ the "part" of $F$ unaffected by this operation. Introducing this basic picture I think about Euclid's "topological" definitions mentioned above. The Flow can be thought of as having an infinite number of dimensions but perhaps it is more appropriate not to prescribe to it any dimension at all. In any event we may prescribe to C *codimension* 1.

Now having the notion of category in mind I represent the obtained construction by the following diagram (Fig.7) leaving the neutral flow $F_C'$ apart:

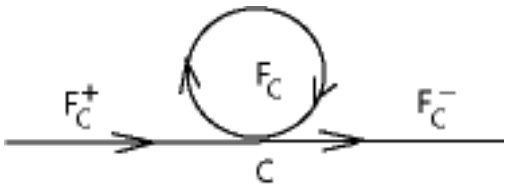

$$F_C^+ \cdot F_C^- = U_C$$

Fig.7

Let us now assume that our basic cut operation is reversible, so that after performing this operation we are still able to "see" how the Flow looked locally before. The kind of reversibility I'm talking about here is usually taken in mathematics for granted. Cutting given sphere $S$ into two semi-spheres $S^+$ and $S^-$, for example, one usually assumes the possibility of bringing the two semi-spheres back together. This allows for thinking of the semi-spheres as two *parts* of the original sphere and writing $S^+ \oplus S^- = S$, where $\oplus$ stands for the mereological sum and $=$ for equivalence (equality), which says that the mereological sum of the two semi-spheres is the same thing as the original sphere. Observe that there is an obvious sense in which $\oplus$ is reverse of cutting, and that if the cutting wouldn't be reversible we wouldn't get anything like equality here.

Similarly writing 2+3=5 we assume that numbers 2 and 3 don't perish forever after being summed up, so 5 can be decomposed into 2 and 3 again. If 2 and 3 would perish after being summed up we could still write 2+3→5 but not 2+3=5. These examples let me think that the kind of reversibility involved here is *necessary* in order to introduce identity (or equality or even equivalence) in mathematics.

As far as we are talking about well-behaved solids and moreover about eternal numbers the reversibility of operations (in the desired sense) is guaranteed but when we try to cope with flows this condition must be mentioned explicitly. Now requiring the reversibility of *cut* we can write $F_C^+ \otimes F_C^- = U_C$ thinking about $\otimes$ as "cut-elimination" and of $U_C$ as a *neighborhood* (local environment) of $C$. The cut-elimination can be equally viewed as *composition* of the incoming and outgoing flows. Noticeably flow $U_C$ resulted from this operation unlike the global flow $F$ is *oriented*, that is, has a certain direction. The property of global flow $F$ that amounts to the possibility of performing the above construction can be thus called *local orientability*. We shouldn't suppose that any given flow is locally orientable. As an example of locally orientable flow think about a *differentiable* flow (say, of incompressible liquid): cuts in this case are spatial points, and their oriented neighborhoods are velocity vectors associated with these points. There are two basic differences between this standard example and our general construction. (i) In a differentiable flow the cut-operation and the cut-elimination are trivialized (and usually not even recognized as operations) through the assumption that points are always "there". In this case what "pours into", "pours out of", and "flows through" any given point is always the same thing, so all local flows are idempotent with respect to their composition. (ii) Unlike points our cuts may have internal dynamics $F_C$, so the whole construction can be reiterated "inside a cut", as we shall shortly see[50].

Notice that operation $\otimes$ just introduced doesn't suggest to regard the equality = as previously given. For it is the reversibility of $\otimes$ and *cut* (each one being the reverse of the other) which allows us to introduce equality here. However we shouldn't suppose that the equality so introduced will be "the" universal mathematical equality. We should rather think about this equality as associated with this particular operation.

---

[50] So the local orientablity is certainly weaker than the differentiability. What about continuity?

Since ⊗ looks very much like the composition of morphisms in a category we can pursue some further categorical-like constructions in $F$. However we can do something more familiar working on the boundary flow $F_C$ which has both "domain" and "codomain". For this end we should think of $F_C$ as a *family* of flows provided with domains and co-domains which are not necessarily identical (remark that the equality associated with ⊗ doesn't apply to $F_C$). Given $F_C$ we perform on it local cuts of the second order (of codimension 2) obtaining usual commutative triangles. Using these triangles we could pursue usual categorical constructions but it seems to be more appropriate to the case to introduce such constructions through cuts of "different shapes" like the square shape. The cuts of the second order can be cut with cuts of the third order, so the construction is continued recursively. What we get looks like n-category constructed downwards. To see this consider a flow of "free" natural transformations, then "cut" these transformation by functors (so the transformations become composable), and finally cut functors with categories (so functors become composable too). Importantly we don't need to stop here and may continue the construction downwards unlimitedly. I cannot develop this idea any further here and leave it for a later study.

## 20. Conclusion

Paraphrasing Quine (1966) we can say: one man's paradox is another man's definition. Considering the invariance through change as a basic feature (if not a definition) of identity we may avoid the *Paradox of Change* but the price will be the lost of primitive and universal character of the identity concept. In my view this price must be paid anyway. People use the word "same" in many different context-dependent senses in everyday talks as well as in scientific discourses. What physicists exactly mean speaking of the "same experiment", "same observation", "same effect, " "same model", "same theory", "same event" or "same particle" ? I don't think that the type/token distinction gives all needed answers. In biology and social sciences things become even more complicated. A philosophical approach to the issue requires first of all distinguishing, specification and theoretic systematization of different senses of the "same" rather than picking one of them, stipulating it as basic and explaining away others. I don't think that Frege is right taking for granted that the notion of identity is unique and simple. I think that Plato was more to the point noticing that nothing like "pure" identity applies to physical and mathematical matters, so he had to stipulate a special realm of eternal Ideas where it might work. But unlike Plato I am rather interested in "impure" identities, which might work in mathematics, physics and other sciences.

We have seen that the issue of identity has been crucial in the development of programs of unification of mathematics since the end of 19th century. Frege's attempts to "fix the identity" of natural numbers as continued by Russell shaped the mainstream philosophy of mathematics (although hardly the mainstream of mathematics) in the XX-th century. The issue of identity remains central in the current program of categorification, which is an alternative project of the unification of mathematics. The fact that our working concept of identity in mathematics is weak and diversified noticed by Plato has been interpreted by Frege, Russell, and their followers as an evidence of the lack of rigor in mathematics, and they tried to fix the problem through introduction of a universal logical notion of identity. Categorification, in contrast, purports to further weakening and diversification of identity revealing genuinely new mathematics in doing so. Interestingly categorification revives certain philosophical ideas, which during the XX-th century remained marginal, like Geach's idea of relative identity and Bradley's "internalism" about relations. It also leads to a repudiation of the idea shared by the majority of Analytic philosophers since Frege that the issue of identity must be firmly fixed from the outset in any serious theoretical enterprise. In the category-theoretic framework the issue of identity is an issue to be studied (both from a general point of view and in every particular case) but not one to be rigidly fixed in advance.

In its present form Category theory doesn't offer yet any systematic account of identity developed in categorical terms. The "internalization of logic" in a category through the standard device of "internal language" doesn't internalize identity. However there are interesting suggestions aiming at a new categorical theory of identity. One is representation of equality of objects in a category as splitting of (categorical) fibration. Another is replacement of equalities in a category by 2-morphisms (morphisms between morphisms) [51]. In both these cases "usual" equalities are still used at certain level of construction. This may be justified through considering this level as "meta-level". However in the case of *n*-categories such approach doesn't look satisfactory. My guess is that the "weakening of identity" wanted in higher categories requires revision of the usual assumption that objects and morphisms in a category form classes. Everybody knows today that the notion of set is not so innocent as it might seem. I tried to show in this paper that the notion of class is not innocent either. I also made an informal suggestion of how to think a category without classes. I believe that Category theory presents an opportunity for philosophers and mathematicians to develop new theories of identity, which will be useful for today's science.

---

[51] A different approach to tackling identity in categories, which I don't discuss here, is presented in (Makkai 1998). Makkai takes a middle way between logicist and mathematist stands and constructs a logical syntactic category "with an appropriate notion of structural equality" as a replacement for the standard equality (p. 179).